\documentclass[11pt]{amsart}
\usepackage{amsmath,amssymb}

\newfont{\Bb}{msbm10 scaled\magstep1}
\newfont{\Bbs}{msbm10 scaled\magstep0}

\newcommand{\rmv}[1]{}

\newcommand{\fq}{\mbox{\Bb F}_{q}}

\newcommand{\fqk}{\mbox{\Bb F}_{q^{k}}}

\newcommand{\sfqk}{\mbox{\Bbs F}_{q^{k}}}
\newcommand{\sfq}{\mbox{\Bbs F}_{q}}

\newcommand{\fp}{\mbox{\Bb F}_{p}}
\newcommand{\fqbar}{\bar{\mbox{\Bb F}}_{q}}
\newcommand{\Tr}{\mbox{\rm Tr}}

\newcommand{\Qp}{\mbox{\Bb Q}_{p}}
\newcommand{\Zp}{\mbox{\Bb Z}_{p}}
\newcommand{\ord}{\mbox{\rm ord}\,}
\newcommand{\sord}{\mbox{\small ord}\,} 
\newcommand{\Oh}{{\mathcal O}}
\newcommand{\SoftOh}{\tilde{{\mathcal O}}} 

\newcommand{\Z}{\mbox{\Bb Z}}
\newcommand{\sZ}{\mbox{\Bbs Z}}
\newcommand{\R}{\mbox{\Bb R}}

\newcommand{\Teich}{\omega} 

\newcommand{\Newt}{\Delta}

\newcommand{\Rram}{R_{1}}
\newcommand{\Runram}{R_{0}}
\newcommand{\Rp}{R}
\newcommand{\extdeg}{a} 
\newcommand{\nvol}{v} 

\newtheorem{theorem}{Theorem}
\newtheorem{corollary}[theorem]{Corollary}
\newtheorem{lemma}[theorem]{Lemma}

\newtheorem{proposition}[theorem]{Proposition}
\newtheorem{definition}[theorem]{Definition}
\newtheorem{algorithm}[theorem]{Algorithm}

\newenvironment{exafont}{\begin{bf}}{\end{bf}}

\renewenvironment{algorithm}{\vspace{0.3cm}\par\noindent\refstepcounter{theorem}\begin{exafont}Algorithm \thetheorem\end{exafont}\hspace{\labelsep}}{\vspace{0.3cm}\par}
\newenvironment{note}{\vspace{0.3cm}\par\noindent\refstepcounter{theorem}\begin{exafont}Note \thetheorem\end{exafont}\hspace{\labelsep}}{\vspace{0.3cm}\par}

\title[Counting points on varieties]
{Counting points on varieties over finite fields
of small characteristic}
\author[Lauder and Wan]{Alan G.B. Lauder and Daqing Wan}

\address{
Computing Laboratory\\
Oxford University\\
Oxford OX1 3QD, UK
\email{}{alan.lauder@comlab.ox.ac.uk}}
\rm{\address{
Department of Mathematics\\
University of California\\ Irvine, CA 92697, USA
\email{}{dwan@math.uci.edu}
}}

\date{\today}

\thanks{
Alan Lauder gratefully acknowledges the support of the EPSRC
(Grant GR/N35366/01) and St John's College, Oxford, and thanks
Richard Brent. Daqing Wan is partially supported by the NSF and
the NSFC. The authors are pleased to thank Colin McDiarmid and
Bernd Sturmfels for answering some questions on convex geometry.
Many excellent suggestions were made for improving the paper by
the anonymous referee, and these were incorporated by the authors
during the revision. They are especially grateful for this help.
}

\begin{document}

\begin{abstract}
We present a deterministic
polynomial time algorithm for computing the zeta function
of an arbitrary variety of fixed dimension over a finite field of small
characteristic. One consequence of this result is an efficient method
for computing the order of the group of rational points on the
Jacobian of a smooth geometrically connected projective
curve over a finite field of small characteristic.
\end{abstract}

\maketitle

\section{Introduction}\label{Sec-Intro}

The purpose of this paper is to give an elementary and self-contained
proof that one may efficiently compute zeta functions of arbitrary varieties
of fixed dimension over finite fields of suitably small characteristic.
This is achieved via the $p$-adic
methods developed by Dwork in his proof of the rationality of the zeta
function of a variety over a finite field \cite{BD60,BD62}.
Dwork's theorem shows that
it is in principle possible to compute the zeta function. Our main
contribution is to show how Dwork's trace formula, Bombieri's
degree bound \cite{EB78} and a semi-linear reduction argument
yield an efficient algorithm for doing so.
That $p$-adic methods may be used to
efficiently compute zeta functions for small characteristic was first
suggested in \cite{DW99,DW00}, where Wan gives a simpler algorithm for
counting the number of solutions to an equation over a finite field
modulo small powers of the characteristic.

We now give more details of our results. For $q = p^{\extdeg}$, where $p$ is
a prime number and $\extdeg$ a positive integer, let $\fq$ denote a finite
field with $q$ elements. Let $\fqbar$ denote an algebraic closure of
$\fq$, and $\fqk$ the subfield of $\fqbar$ of order $q^{k}$.
Denote by $\fq [X_{1},\dots,X_{n}]$ the
ring of all polynomials in $n$ variables over $\fq$.

For a polynomial $f \in \fq[X_{1},\dots,X_{n}]$, we
denote by $N_{k}$ the number of solutions to the equation
$f = 0$ with coordinates in $\fqk$. The zeta function
of the variety defined by $f$
is the formal power series in $T$ with non-negative integer coefficients
\[ Z(f/\fq)(T) = \exp \left(\sum_{k = 1}^{\infty}
 \frac{N_{k}T^{k}}{k}\right).\]
Dwork's theorem asserts that $Z(f/\fq)$ is a rational function
$r(T)/s(T)$ with integer coefficients. From this it follows that
knowledge of explicit bounds $\deg(r) \leq D_{1}$ and $\deg(s)
\leq D_{2}$, and of the values $N_{k}$ for $k = 1,2,\dots,D_{1} +
D_{2}$ is enough to efficiently determine $Z(f/\fq)$, see
\cite{DW00}. The Bombieri degree bound tells us that $\deg(r) +
\deg(s) \leq (4d + 9)^{n+1}$, see \cite{EB78,DW00}, and so in
particular we may take $D_{1} = D_{2} = (4d+9)^{n+1}$. Each number
$N_{k}$ can be computed in a naive fashion by straightforward
counting, using $q^{nk}$ evaluations of the polynomial $f$. Thus
one may compute the zeta function of a variety, but the naive
method described requires a number of steps that is exponential in
the parameters $d^{n}$ and $\log q$, where $d$ is the total degree
of $f$. Note that the dense input size of $f$ is $\Oh ((d+1)^n\log
q)$, and the size of the zeta function is polynomial in
$\Oh((d+1)^n \log q)$, by the Bombieri degree bound. We prove the
following theorem.

\begin{theorem}\label{Thm-Main}
There exist an explicit deterministic algorithm and an explicit polynomial $P$
such that for any $f \in \fq[X_{1},\dots,X_{n}]$ of total degree
$d$, where $q = p^{\extdeg}$ and $p$ is prime, the algorithm computes the
zeta function $Z(f/\fq)(T)$ of $f$ in a number of bit operations
which is bounded by $P(p^{n}d^{n^{2}}\extdeg^{n})$.
\end{theorem}

In particular, this computes the zeta function of a polynomial
in a fixed number of variables over a finite field of ``small
characteristic'' in deterministic polynomial time.
Note that our result makes no assumption of
non-singularity on the variety defined by the polynomial. Also, we
shall explicitly describe all the algorithms in this paper, rather
than just prove their existence, and we assume that the finite field
$\fq$ itself is presented as input via an irreducible polynomial of
degree $\extdeg$
over the prime field $\fp$, as explained in Section \ref{Sec-PadicTheory}.

With regard to the exponents in the algorithm, we state these
precisely in Theorem \ref{Thm-MainRefined}.  For now we observe
that if $p$, $d$ and $n$ are fixed then the time required to
compute the number of points in $\fq$ is $\Oh(\extdeg^{3n + 7})$,
with space complexity $\Oh(\extdeg^{2n+4})$. Here we are ignoring
logarithmic factors (see also Proposition \ref{Prop-BitCompTotal}).
Note that all our complexity estimates are made using standard
methods for multiplication in various rings, and can be modestly
reduced with faster methods.

We also present refinements to this result based upon the ideas of
Adolphson and Sperber \cite{AS87}, and indeed from the outset will
follow their approach, as it involves little extra complication.
This refinement takes into account the terms which actually occur in the
polynomial $f$ rather than working solely with the total degree.
We shall need more definitions: the
support of a polynomial $f$ is the set of exponents $r =
(r_{1},\dots,r_{n})$ of non-zero terms $a_{r}X_{1}^{r_{1}}\dots
X_{n}^{r_{n}}$ which occur in $f$, thought of as points in $\R^{n}$.
The Newton polytope of $f$
is defined to be
the convex hull in $\R^{n}$ of the support of $f$.
Our refined version of Theorem \ref{Thm-Main} essentially replaces the
parameter $d^{n}$ with the normalised volume of the Newton
polytope of $f$ (see Proposition \ref{Prop-BitComp} and Section
\ref{Sec-BitCompAux}).

Zeta functions may also be defined for a finite collection of polynomials,
and we next describe
how our results may be extended to this case. An affine variety
$V$ over $\fq$ is the set of common zeros in $(\fqbar)^{n}$ of a
set of polynomials $f_{1},f_{2},\dots,f_{r}$. An analogous
zeta function $Z(V/\fq)$ may then be defined in terms of the number
of solutions in each finite extension field of $\fq$. These numbers may be
computed using an inclusion-exclusion argument involving the
polynomials $\prod_{i \in S} f_{i}$ where $S$ is a subset of $\{1,2,\dots
,r\}$, see \cite{DW00}. Theorem \ref{Thm-Main} then easily yields the
following.

\begin{corollary}\label{Cor-GeneralVar}
There exist an explicit deterministic algorithm and an explicit
polynomial $Q$ with the following property. Let $f_{1},\dots,f_{r}
\in \fq[X_{1},\dots,X_{n}]$ have total degrees $d_{1},\dots,d_{r}$
respectively, where $q = p^{\extdeg}$ and $p$ is prime. Denote by $V$
the affine variety defined by the common vanishing of these
polynomials, and define $d = \sum_{i = 1}^{r}d_{i}$.  The algorithm
computes the zeta function $Z(V/\fq)(T)$ in a number of bit operations
which is bounded by $Q(p^{n}d^{n^{2}}\extdeg^{n} 2^{r})$.
\end{corollary}

Thus one has an efficient algorithm for computing the zeta function of
an arbitrary affine variety over $\fq$ assuming the characteristic, dimension and the number
of defining polynomials are fixed. More generally still, an arbitrary
variety is defined through patching together suitable affine
varieties. Zeta functions for such general varieties may be
defined. These zeta functions may be computed using the above ideas
provided explicit data are given on how to construct them from affine
patches. As an example, the zeta functions of arbitrary projective
varieties or toric varieties may be computed in this way.

Conceptually our algorithm is rather straightforward. The zeta
function can be expressed in terms of the ``characteristic
power series'' (Fredholm determinant) of a
certain ``lifting of Frobenius'' which acts on an infinite
dimensional $p$-adic Banach space constructed by Dwork. Under modular
reduction, we obtain an operator acting on a finite dimensional
vector space. Unfortunately, this operator cannot be computed efficiently
directly from its definition;
however, it can be expressed as a product of certain
semi-linear operators each of which can be computed efficiently if the
characteristic $p$ is small.  This last step is thus of crucial
importance in deriving
an efficient algorithm. Note that the same idea is used in
\cite{DW99,DW00} in a simpler situation.  In more concrete language,
the algorithm requires one to construct a certain ``semi-linear''
finite matrix and compute the ``linear'' matrix which is the product
of the Galois conjugates of the semi-linear matrix.  The number of
points is then read off from the trace of the final linear matrix. The
zeta function can be computed from the characteristic polynomial of
the final linear matrix.
The semi-linear matrix itself is defined
over a certain finite ``$p$-adic lifting'' of the original finite
field.

In the literature algorithms have already been developed for computing
zeta functions of curves and abelian varieties \cite{AH96, NE98, JP90, RS85,
RS95}. These utilise the theory originally developed by Weil for
abelian varieties, whereas we use Dwork's more general and simpler
$p$-adic theory.
For example, in the case of a smooth geometrically irreducible projective plane
curve of degree $d$ over a field of
size $q=p^{\extdeg}$ these algorithms have a time complexity which grows as
$(\log{q})^{C_d}$, where the constant $C_d$
grows exponentially in the degree $d$.
Given an absolutely irreducible bivariate polynomial $f$ of degree $d$
over $\fq$,
the zeta function of the unique smooth projective curve birational
to the affine curve defined by $f$ may be computed in time polynomial
in $d$, $p$ and $\extdeg$ using our approach.
Thus our more general method
is far better in terms of the degree $d$ if the characteristic $p$ is
small, since the running time has polynomial growth in $d$ (but much
worse if $p$ is large).
The zeta function immediately gives the order of the group of rational
points on the Jacobian,
see \cite{DW00}.

\begin{corollary}\label{Cor-Jacobian}
There exist an explicit deterministic algorithm and an explicit
polynomial $R$ with the following property.  Let $V$ be a
geometrically irreducible affine curve defined by the vanishing of
polynomials $f_1,\dots,f_r \in \fq[X_{1},\dots,X_{n}]$ of total
degrees $d_{1},\dots,d_{r}$ respectively,
where $q = p^{\extdeg}$ and $p$ is prime.  Denote by $\tilde{V}$ the
unique smooth projective curve birational to the affine curve
$V$, and let $d = \sum_{i = 1}^{r}d_{i}$.
The algorithm computes the order of the group of rational
points on the Jacobian of $\tilde{V}$ in a number of bit operations
bounded by $R(p^nd^{n^2}\extdeg^n 2^r)$.
\end{corollary}

Thus one may compute the order of the group of rational points on
the Jacobian of a smooth geometrically irreducible projective
curve over a finite field of small characteristic
in deterministic polynomial time,
provided
the number of variables and the number of defining equations are fixed.
(Notice that the case $r = 1$ and $n = 2$ corresponds to that
of being given a possibly singular plane model of the curve.)
In particular, this answers a question posed in \cite{BP96}, which Poonen
attributes to Katz and Sarnak. We note that recently a similar
result for special classes of plane curves was independently obtained
in \cite{GG01,KK01}, using Monsky-Washnitzer's method.
Also, a different $p$-adic approach for elliptic curves has been
developed in \cite{TS00}.

This paper is written primarily for theoretical computational
interest, in obtaining a deterministic polynomial time algorithm for
computing the zeta function in full generality if $p$ is small.  It
can certainly be improved in many ways for practical computations.
What we have done in this paper is to work on the easier but more flexible
``chain level''. A general improvement in the smooth case,
is to work on the cohomology level. Then there are several related
$p$-adic cohomology theories available, each
leading to a somewhat different version of the algorithm.  Again, as
indicated in \cite{DW99,DW00}, these $p$-adic methods are expected to
be practical only for small $p$. (The authors have
recently applied Dwork's cohomology theory to derive a practical
algorithm in a special case \cite{LWLMS}.)

\section{Additive character sums over finite fields}

The most natural objects of study in Dwork's theory are certain
additive character sums over finite fields. In this section we introduce
these sums, and explain their connection to varieties.

An additive character $\Psi$ is a mapping from $\fqk$ to the
group of units of some commutative ring $S$ with identity $1$ such
that $\Psi(x + y) = \Psi(x)\Psi(y)$ for $x,y \in \fqk$. We say
that it is non-trivial if $\Psi(x) \ne 1$ for some $x \in \fqk$.
Let $\{\Psi_{k}\}_{k \geq 1}$ be any family of mappings with each
$\Psi_{k}$ a non-trivial additive character from $\fqk$ to some
extension ring of the integers whose image is a group of order $p$
with elements summing to zero. We assume that the family $\{
\Psi_k\}$ forms a tower of characters in the sense that for each
$k\geq 2$,
$$\Psi_k =\Psi_{1}\circ {\rm Tr}_{\sfqk/\sfq}.$$
In our application, the ring $S$ will be
taken to be a certain $p$-adic ring.

For the remainder of the paper we will
use multi-index notation.
Specifically, we
let $X^{u}$ represent the monomial
$X_{0}^{u_{0}}X_{1}^{u_{1}}\dots X_{n}^{u_{n}}$ for an integer vector
$u = (u_{0},u_{1},\dots,u_{n})$; let $x$ be an $(n+1)$-tuple
$(x_{0},x_{1},\dots,x_{n})$ of field elements; and
$X$ the list of indeterminates $X_{0},X_{1},\dots,X_{n}$.
Observe here that we have introduced an extra indeterminate
$X_{0}$. Even though the polynomial $f$ whose zeta function
we wish to compute is in the $n$ variables $X_{1},\dots,X_{n}$,
in Dwork's theory the extra indeterminate arises naturally, as we
are about to see.

\begin{lemma}\label{Lem-ExpSum}
Let $f \in \fq[X_{1},\dots,X_{n}]$ and $N_{k}^{*}$ denote the
number of solutions to the equation $f = 0$ in the affine
torus $(\fqk^{*})^{n}$. Then
\[ \sum_{x \in (\sfqk^{*})^{n+1}} \Psi_{k}(x_{0}f(x_{1},\dots,x_{n}))
 = q^{k}N^{*}_{k} - (q^{k} - 1)^{n}\]
where $x = (x_{0},x_{1},\dots,x_{n})$.
\end{lemma}

\begin{proof}
We have that for any $u \in \fqk$
\[ \sum_{x_{0} \in \sfqk} \Psi_{k}(x_{0}u) =
\left\{
\begin{array}{ll}
0 & \mbox{if $u \in \fqk^{*}$}\\
q^{k} & \mbox{if $u = 0$}.
\end{array}
\right.
\]
This is a standard result from the theory of additive character sums
\cite[page 168]{LN86}.
Thus $\sum \Psi_{k}(x_{0}f(x_{1},\dots,x_{n}))$, where
the sum is taken over points $x \in \fqk \times (\fqk^{*})^{n}$,
equals $q^{k}N_{k}^{*}$. Removing the contribution of $(q^{k} - 1)^{n}$
from the terms with $x_{0} = 0$ in this sum gives the required result.
\end{proof}

Our next step will be to find an alternative formula for the
lefthand side of the equation in Lemma \ref{Lem-ExpSum}. This is
achieved in Proposition \ref{Prop-AnalyticExpression}, which leads us
eventually to Dwork's Trace Formula (Theorem \ref{Thm-DTF}).

\section{p-Adic Theory}\label{Sec-PadicTheory}

\subsection{p-adic rings}
We first introduce notation for the $p$-adic rings we shall need,
before explaining how to construct them and compute in them. Let
$\Qp$ be the field of $p$-adic rationals, and $\Zp$ the ring of
$p$-adic integers (see \cite{NK}). Denote by $\Omega$ the
completion of an algebraic closure of $\Qp$. Select $\pi \in
\Omega$ with $\pi^{p-1} = -p$ and define $\Rram = \Zp[\pi]$, a
totally ramified extension of $\Zp$ of degree $p-1$. By binomial
expansion and Hensel's lifting lemma, one sees that the equation
$(1+\pi t)^p =1$ has exactly $p$ distinct solutions $t$ in
$\Rram$. In particular, $\Rram$ contains all $p$-th roots of
unity. The motivation behind the introduction of $\Rram$ is that
to define an additive character of order $p$, we need a small
$p$-adic ring which contains a primitive $p$th root of unity.

Let $\Runram$ denote the ring of integers
of the unique unramified extension of
$\Qp$ in $\Omega$ of degree $\extdeg$, where $q = p^{\extdeg}$. Finally, let
$\Rp$ be the compositum ring of $\Rram$ and $\Runram$. We have the
diagram of ring extensions

\[
\begin{array}{ccccc}
& & \Omega & &\\
& & | & & \\
 & & \Rp & &\\
& \diagup & & \diagdown &\\
\Runram & & & & \Rram \\
& \diagdown & & \diagup &\\
& & \Zp & &
\end{array}
\]
The residue class ring of $\Runram$ is $\fq$, and we shall
``lift'' the coefficients of the polynomial $f \in
\fq[X_{1},\dots,X_{n}]$ to this characteristic zero ring, using
the Teichm\"uller lifting. Precisely, the Teichm\"{u}ller lift
$\Teich(x)$ of a non-zero element $x \in \fqbar$ is defined as the
unique root of unity in the maximal unramified extension of $\Qp$
which is congruent to $x$ modulo $p$ and has order coprime to $p$.
We define $\Teich(0) = 0$. Thus the compositum $\Rp$ contains both
the lifting of the coefficients of $f$ and the image of an
additive character we shall construct.

\subsection{Algorithmic aspects}

\subsubsection{Construction and lifting Frobenius}\label{Sec-ContLiftFrob}

We assume that $\fq$ is presented as the quotient $\fp[y]/(h)$ where $h(y)$
is a monic, irreducible polynomial of degree $\extdeg$
over the prime field $\fp$.  For
any positive integer $N$ we describe how the quotient ring
$\Rp/(p^{N})$ may be constructed:
Lift the polynomial $h$ to an integer polynomial
$\hat{h}$ whose coefficients lie in the open interval $(-p/2, (p+1)/2)$.
Take the unramified extension $\Runram$ to be $\Zp[\mu]=\Zp[y]/(\hat{h})$.
Elements in $\Runram/(p^{N})$ can now be represented as linear combinations
over $\Zp/(p^{N})$ of the basis elements
$1,\mu,\dots,\mu^{\extdeg-1}$.
(Recall that $\Zp/(p^{N})$ can be identified with $\Z/(p^{N})$.)
The extension
$\Rp/(p^{N}) = \Runram[\pi]/(p^{N})$ is easily constructed
by adjoining an element $\pi$ and specifying the
relation $\pi^{p-1} = -p$.

We define a lifting of the Frobenius automorphism on $\fq$
to an automorphism of $\Rp$ which is the identity on $\Rram$. Define
the map $\tau:\Rp \rightarrow \Rp$
by setting $\tau(\mu)$ to be the unique root
of $\hat{h}$ which is congruent to $\mu^{p}$ modulo $p$. Define
$\tau(\pi) = \pi$, and extend to the whole of $\Rp$ by
insisting $\tau$ is an automorphism.
We extend $\tau$ to act on
$\Rp[[X]]$ coefficient-wise, fixing monomials.
Here $\Rp[[X]]$ is the ring of all power series
in the indeterminates $X = X_{0},\dots,X_{n}$ with coefficients
from $\Rp$.

\subsubsection{Complexity of arithmetic}
In this section we bound the complexity of the basic arithmetic
operations in the ring $\Rp/(p^{N})$, along with that of computing the
map $\tau$ and Teichm\"{u}ller lifts. Note that these
estimates are all the simplest possible, and can be improved
using more advanced methods. The reader may wish to skip the
proof of the next lemma, and refer back when required in Section
\ref{Sec-BitCompRing}.

\begin{lemma}\label{Lem-R0ops}
Elements in
$\Rp/(p^{N})$ can be represented using
$\Oh(p\extdeg N\log{p})$ bits. Addition and subtraction
can be performed in
$\Oh(p\extdeg N\log{p})$
bit operations, and
multiplication and inversion of units in
$\Oh((p\extdeg N\log{p})^{2})$ bit operations.
The Teichm\"{u}ller lifting to
$\Rp/(p^{N})$ of a finite field element can be computed
in $\Oh((\extdeg\log{p})^{3}N^{2})$ bit operations.
For any $1 \leq i \leq \extdeg - 1$,
the map $\tau^{i}$ on $\Rp/(p^{N})$ may be evaluated using
$\Oh(p(\extdeg N\log{p})^{2})$ bit operations. (For
the powers of $\tau$ we require a total of
$\Oh(\extdeg^{4}N^{2}(\log{p})^{3})$ bits of precomputation.)
\end{lemma}

\begin{proof}
Elements in $\Rp/(p^{N})$ can be written as
\begin{equation}\label{Eqn-RingElt}
\sum_{i = 0}^{\extdeg -1}\sum_{j = 0}^{p-2} c_{ij}\pi^{j}\mu^{i}
\end{equation}
where the coefficients $c_{ij}$ belong to the ring
$\Zp/(p^{N})$ of size $p^{N}$. The bit size of such an expression
is $\extdeg(p-1)\log(p^{N}) = \Oh(p\extdeg N\log{p})$.
Addition and subtraction
are straightforward, just involving the addition of integers and
reduction modulo $p^{N}$. Likewise, multiplication of two expansions
of the form (\ref{Eqn-RingElt}) is straightforward, using the
reduction relations $\pi^{p-1} = -p$ and $\hat{h}(\mu) = 0$.

For Teichm\"{u}ller lifting and inversion we shall use Newton
iteration with quadratic convergence.
Specifically, we shall use Newton lifting with respect to the prime
$\pi$ in $\Rp$, and define $l=(p-1)N$ so that $\pi^{l} = (-p)^{N}$.
Given a polynomial $\phi(Y) \in
\Rp[Y]$ and an element $g_{0} \in \Rp$ such that $\phi(g_{0})
 \equiv 0 \bmod{\pi}$ and $\phi^{\prime}(g_{0})$ is invertible
(modulo $\pi$) with inverse $s_{0}$ modulo $\pi$,
this algorithm
computes an element $g \in \Rp/(\pi^{l})$ such that $\phi(g) \equiv
0 \bmod{\pi^{l}}$ and
$g \equiv g_0 \bmod{\pi}$ (compare with \cite[Algorithm 9.22]{GaGe}).
For $i \geq 1$, assuming that $g_{i-1}$ and $s_{i-1}$ have been found,
we define
$g_{i} = g_{i-1} - \phi(g_{i-1})s_{i-1} \bmod{\pi^{2^{i}}}$ and
$s_{i} = 2s_{i-1} - \phi^{\prime}(g_{i})s_{i-1}^{2} \bmod{\pi^{2^{i}}}$.
As in \cite[Theorem 9.23]{GaGe} one checks that
$g_{i} \equiv g_{0} \bmod{\pi}$, $\phi(g_{i}) \equiv 0
\bmod{\pi^{2^{i}}}$ and $s_{i} \equiv \phi^{\prime}(g_{i})^{-1}
\bmod{\pi^{2^{i}}}$ at the $i$th step. Thus after $\lceil \log_{2}(l)
\rceil$ steps we shall have found the required approximate root $g =
g_{\lceil \log_{2}(l) \rceil}$. The $i$th step involves three
additions/subtractions
and three multiplications in the ring $\Rp/(\pi^{2^{i}})$, which
requires $\Oh((2^{i}\extdeg \log{p})^{2})$ bit operations, along
with evaluation of the polynomials $\phi$ and $\phi^{\prime}$ modulo
$\pi^{2^{i}}$. Let $c(\phi,i)$ be the complexity of these latter operations.
Thus the
total complexity is $\Oh((p\extdeg N\log{p})^{2} + c(\phi))$ bit operations
where $c(\phi) = \sum_{i = 1}^{\lceil \log_{2}(l) \rceil} c(\phi,i)$.
Note that to compute approximate roots in the ring $\Runram$ rather
than $\Rp$, one can lift using the prime $p$ rather than $\pi$ and
obtain a complexity of $\Oh((\extdeg N\log{p})^{2} + c(\phi))$.

Suppose now we are given a unit $x \in \Rp/(p^{N})$. We compute
a Newton lifting
starting from an element $g_{0} \in \Rp/(\pi)$ such that $xg_{0}
\equiv 1 \bmod{\pi}$; that is, we use the equation $\phi(Y) = xY -
1$. Now $\Rp/(\pi)$ is just the finite field $\fq$ and so an inverse
$g_{0}$ of $x$ modulo $\pi$ can be computed in $\Oh((\extdeg \log{p})^{2})$
bit operations \cite[Corollary 4.6]{GaGe}. (Here
$\phi^{\prime}(Y) = x$ and $s_{0} = g_{0}$, and so in fact we only
need to iterate the formula for $s_{i}$.)
In this case
$c(\phi,i)$ is just one multiplication  and
a subtraction modulo $\pi^{2^{i}}$.
By the above paragraph the
total complexity is then $\Oh((p\extdeg N\log{p})^{2})$ for inversion.
For Teichm\"{u}ller lifts we use the same approach, only with the
polynomial $\phi(Y) = Y^{q-1}  - 1$ and lifting in $\Runram$ via the element
$p$. Here using a fast exponentiation routine we find that
$c(\phi,i)$ involves $\Oh(\log{q})$ multiplications and a subtraction in
$\Runram/(p^{2^{i}})$. Thus $c(\phi) = \Oh((\extdeg N\log{p})^{2}\log{q})$
which gives the Teichm\"{u}ller lifting estimate.

We precompute $\tau$ on the basis elements $1,\mu,\dots,
\mu^{\extdeg -1}$.
To do this, recall that $\tau(\mu) \in \Runram$ is defined as the
unique root of $\hat{h}$ which is congruent modulo $p$ to
$\mu^{p}$. This may be approximated modulo $p^{N}$ by Newton lifting
in $\Runram$ with respect to $p$ using the polynomial $\phi(Y) =
\hat{h}(Y)$ and the initial value $g_{0} = \mu^{p}$.
(Finding $\mu^{p}$ takes $\Oh(a^{2}(\log{p})^{3})$ bit operations
and this is absorbed in the stated precomputation estimate.)
Here
$c(\phi,i)$ is $\Oh(\extdeg)$ additions and multiplications in
$\Runram/(p^{2^{i}})$, using Horner's method for polynomial evaluation
\cite[Page 93]{GaGe}. Thus the Newton lifting estimate gives a
complexity of $\Oh(\extdeg(\extdeg N\log{p})^{2})$ bit operations. Using
$\tau(\mu^{i}) = (\tau(\mu))^{i}$ one can now find the image of all
basis elements in a further $\Oh(\extdeg(\extdeg
N\log{p})^{2})$ bit operations.
One stores this information as a matrix for $\tau$ acting
on $\Runram/(p^{N})$ as an $\Zp/(p^{N})$-module with basis the powers
of $\mu$. The map $\tau$ can now be computed on any element in
$\Runram/(p^{N})$ in $\Oh(\extdeg^{2}(N\log{p})^{2})$ bit operations using
linear algebra. By taking powers of the matrix, matrices for the
maps $\tau^{i}$ for $1 \leq i \leq \extdeg - 1$ can also be found
in $\Oh(\extdeg^{3}\extdeg(N\log{p})^{2})$ bit operations.
Thus the total precomputation is bounded by $\Oh(\extdeg^{4}N^{2}
(\log{p})^{3})$, and each evaluation of $\tau^{i}$ on
$\Runram/(p^{N})$ takes $\Oh((aN\log{p})^{2})$ bit operations.
Finally, to evaluate $\tau$ on $\Rp/(p^{N})$ one writes elements of
$\Rp$ on the $\Runram$-basis $1,\pi,\dots,\pi^{p-2}$ and applies
$\tau$ component-wise.
\end{proof}

\subsection{$p$-adic valuations and convergence of power series}

Denote by $\ord$ the additive valuation on $\Omega$ normalised so
that $\ord(p) = 1$. Thus $\ord(\pi) = 1/(p-1)$.
Define a $p$-adic norm $|.|_{p}$ on $\Omega$ by
$|x|_{p} = p^{-\sord(x)}$. The set of all $x \in \Omega$ with
$|x|_{p} \leq 1$ (equivalently $\ord(x) \geq 0$) is called the
closed unit disk.
Given any formal power series $\sum_{r} A_{r}
X^{r}$ where $A_{r} \in \Omega$, we say it converges at a point $x
= (x_{0},\dots,x_{n}) \in \Omega^{n+1}$ if the sequence of
partial sums $\sum_{r,\,|r| < e} A_{r} x^{r}$ tends to a limit under
the $p$-adic norm. (Here $|r| = \sum_{i = 0}^{n} r_{i}$.) This
sequence of partial sums will converge if and only if the summands
$A_{r}x^{r}$ tend to zero $p$-adically (that is, are divisible in
the ring of integers of $\Omega$
by increasingly large powers of $p$), as $|r|$ goes to
infinity. In particular, if there is a real number $c > 0$ such that
$\ord(A_{r}) \geq c|r|$ for all
$r$, then the series will certainly converge for all points $x$
which are Teichm\"{u}ller liftings of points over $\fqbar$.
Note that throughout the paper we shall use the
additive valuation $\ord$ rather than the $p$-adic norm $|.|_{p}$
itself.

\section{Analytic representation of characters}\label{Sec-ARC}

\subsection{Dwork's splitting functions}

We now need to find a suitable $p$-adic expression for a non-trivial
additive character from $\fq$ to $\Rram$.  In the case of complex
characters, this is done via the exponential function.
However, the radius of convergence of the exponential
function in $\Omega$ is not large enough; in particular,
it does not converge on the Teichm\"{u}ller lifting of all the
points in $\fq$. Instead we use the power series constructed by Dwork
using the exponential function. (The reader may find the
discussion of the related Artin-Hasse function on \cite[Pages
92-93]{NK} helpful.)

Let $\mu$ be the M\"obius function.
Taking the logarithmic derivative, one checks that the exponential function
has the following product expansion
\[
\exp(z)=\sum_{k=0}^{\infty}\frac{z^k}{k!} = \prod_{k=1}^{\infty} (1-z^k)^{-\mu(k)/k}.\]
This can be rewritten as
\[
\exp(z)=\prod_{(k,p)=1}^{\infty} (1-z^k)^{-\mu(k)/k}(1-z^{kp})^{\mu(k)/(kp)}.
\]
It follows that
\begin{equation}\label{Eqn-ProdExp1}
\exp\left(z +\frac{z^p}{p}\right) = \prod_{(k,p)=1}^{\infty}
(1-z^k)^{-\mu(k)/k}(1-z^{kp^2})^{\mu(k)/(kp^2)}.
\end{equation}
Replacing $z$ by $\pi z$ in the above
relation and noting that $\pi^p =-p\pi$, we define a power series in $z$
by
\[\theta(z) = \exp(\pi z - \pi z^p).
\]
Writing $\theta(z) = \sum_{r = 0}^{\infty}\lambda_{r}z^{r}$ we see that
$\lambda_{r} = \pi^{r}/r!$ for $r < p$, and we shall shortly show
that all $\lambda_{r}$ lie in $\Rram$.
From (\ref{Eqn-ProdExp1}) we get the product expansion
\begin{equation}\label{Eqn-ProdExp2}
\theta(z)=\prod_{(k,p)=1}^{\infty}
(1-\pi^k z^k)^{-\mu(k)/k}(1-\pi^{kp^2}z^{kp^2})^{\mu(k)/(kp^2)}.
\end{equation}
By the binomial expansion, the first factor
\[
(1-\pi^kz^k)^{-\mu(k)/k}=\sum_{j=0}^{\infty}(-1)^j{{-\mu(k)/k}\choose j}\pi^{jk}z^{jk}
=\sum_{j=0}^{\infty}b_j(k)z^{kj}
\]
is a power series in $z^k$.
Now for $(k,p)=1$, we have \cite[Page 82]{NK}
\[
{-\mu(k)/k}\in \Z_p, \ \ {{-\mu(k)/k}\choose j} \in \Z_p.
\]
Thus for $j > 0$ the coefficient of $z^{jk}$ satisfies
\[
\ord(b_j(k)) \geq \frac{jk}{p-1} >  \frac{p-1}{p^2} jk.
\]
Similarly, for the second factor, we write
\[
(1-\pi^{kp^2}z^{kp^2})^{\mu(k)/(kp^2)}
=\sum_{j=0}^{\infty}(-1)^j{{\mu(k)/kp^2}\choose j}(\pi z)^{jkp^2}
=\sum_{j=0}^{\infty} c_j(kp^2)z^{jkp^2}.
\]
Now for $j > 0$ we have $\ord(j!) <
\frac{j}{p-1}$, see \cite[Page 79]{NK}. Thus
for $j > 0$ and $k$ positive and coprime to $p$,
the coefficient of
$z^{jkp^2}$ satisfies
\[
\ord(c_j(kp^2)) > \frac{jkp^2}{p-1} -2j -\frac{j}{p-1} \geq  \frac{p-1}{p^2} jkp^2.
\]
Putting the above two inequalities together, we conclude that
for $r > 0$ the coefficients $\lambda_{r}$ of $\theta(z)$ satisfy
\begin{equation}\label{Eqn-DecayLambdas}
\ord(\lambda_{r}) > \frac{(p-1)r}{p^{2}},\,\, \lambda_{r} \in \Rram.
\end{equation}
This shows that the power series $\theta(z)$ is convergent in the disk
$|z|_p<1+\epsilon$ for some $\epsilon >0$. In particular, $\theta(z)$
converges on the closed unit disk, and Definition \ref{Def-theta}
makes sense.

In the proof of the next lemma we shall use the fact that
$\ord(\lambda_{r}) \geq 2/(p-1)$ for $r \geq 2$, and so
\begin{equation}\label{Eqn-theta1modpi2}
\theta(z) \equiv 1 + (\pi z) \bmod{(\pi z)^{2}}.
\end{equation}
This can be seen as follows: The Artin-Hasse exponential
function \cite[Page 93]{NK}
\[ E(z):= \prod_{(k,p) = 1}^{\infty}(1 - z^{k})^{-\mu(k)/k} =
\exp\left(z + \frac{z^{p}}{p} + \frac{z^{p^{2}}}{p^{2}} + \dots
\right)\]
has coefficients in $\Zp$, because each factor in the product expansion
does. Since $E(\pi z) \equiv \theta(z) \bmod{z^{p^{2}}}$ we see
\[ \ord(\lambda_{r}) \geq \frac{r}{p-1}\]
for $0 \leq r < p^{2}$. This estimate combined with
(\ref{Eqn-DecayLambdas}) gives (\ref{Eqn-theta1modpi2}).

\begin{definition}\label{Def-theta}[Dwork's splitting function]
Let
\[
\Phi_{k}(z) = \prod_{i = 0}^{\extdeg k-1}\theta(z^{p^{i}}) \in
\Rram[[z]],\]
and
\[
\Psi_{k} = \Phi_{k} \circ \Teich: \fqk \rightarrow
\Rram, \]
where $\Teich$ is the
Teichm\"{u}ller map. (Recall that $q = p^{a}$.)
\end{definition}

(That $\Psi_k$ has image in $\Rram$ can be seen as follows:
Let $\Rram[\Teich(\fqk)]$ be $\Rram$ adjoined the image of
$\Teich$ on $\fqk$. Then $\Rram[\Teich(\fqk)]$ is an unramified
extension of $\Rram$ of degree $k$. The Galois group of the
corresponding quotient field extension is generated by
$\tau$. The map $\tau$ acts on a Teichm\"{u}ller point $\Teich(x)$
as $\tau(\Teich(x)) = \Teich(x)^p$. Hence it fixes the element
$\Psi_k(x)$ for $x \in \fqk$, and so $\Psi_k(x) \in \Rram$.)

\begin{lemma}\label{Lem-DSF}
The maps $\Psi_{k}$ form a tower of non-trivial additive
characters from the fields $\fqk$ to the ring $\Rram$.
\end{lemma}

\begin{proof}
(This is the case ``$s=1$'' on \cite[Pages 55-57]{BD62}.)
We first show that $\theta(1)$ is a primitive $p$th root of unity.
By (\ref{Eqn-theta1modpi2}) we see $\theta(1) \ne 1$.
As a formal power series in $z$,
\[
\theta(z)^p = \exp(p\pi z)\exp(-p\pi z^p).
\]
Now, $\theta(z)$, $\exp(p\pi z)$ and $\exp(-p \pi z^p)$
are all convergent in $|z|_{p} < 1 + \varepsilon$ for some
$\varepsilon > 0$. We can thus substitute $z=1$ and find that
$\theta(1)^{p} = \exp(p\pi) \exp(-p \pi)= 1$.
Thus $\theta(1)$ is a primitive $p$-th root of unity in $\Rram$.

Next, for
$\gamma \in \Runram$ with $\gamma^{p^{\extdeg k}} = \gamma$, we claim that
\[ \prod_{i = 0}^{\extdeg k - 1}\theta(\gamma^{p^{i}}) =
\theta(1)^{\gamma + \gamma^{p} + \dots + \gamma^{p^{\extdeg k-1}}}.\]
Using (\ref{Eqn-theta1modpi2}) it is clear that both sides are congruent to
\[
1+\pi (\gamma + \gamma^{p} + \dots + \gamma^{p^{\extdeg k-1}})
\]
modulo $\pi^2$. To prove the claim, it
remains to prove that both sides are $p$-th
roots of unity. The right side is a $p$-th root of unity since
$\theta(1)$ is a $p$-th root of unity. The $p$-th power of the
left side is
\[
\prod_{i = 0}^{\extdeg k - 1}\theta(\gamma^{p^{i}})^p
=\exp(p \pi \sum_{i=0}^{\extdeg k-1}(\gamma^{p^{i}} - \gamma^{p^{i+1}}))
=\exp(p \pi \gamma)\exp(- p \pi \gamma^{p^{\extdeg k}})=1.
\]
Thus, the left side is also a $p$-th root of unity. The claim is proved.
Note that the individual factor $\theta(\gamma)$ is not necessarily
a $p$-th root of unity.
We conclude that
\[ \Psi_{k}(x) = \theta(1)^{\mbox{\small Tr}_{k}(x)}\]
for any $x \in \fqk$, where $\Tr_{k}$ is the
trace function from $\fqk$ to $\fp$, and
the exponent is thought of as an integer.
\end{proof}

\begin{note}
The infinite sum
\[
\exp(\pi (z-z^p))=\sum_{k=0}^{\infty}\frac{(\pi(z- z^p))^k}{k!}
\]
is convergent for $|z|_p<1$, but not necessarily convergent
for $|z|_p=1$. If $|z|_p=1$, then it is possible that $|z-z^p|_p=1$
and for such $z$ the above infinite sum does not converge.
Since the infinite sum does not converge everywhere on the
disk $|z|_{p} \leq 1$, one cannot simply
substitute $z=1$ into the
above infinite sum and get the contradiction that $\theta(1)=1$.
There is no contradiction here!
\end{note}

We now define another power series related to our original
polynomial $f$ whose relevance will become apparent in
Proposition \ref{Prop-AnalyticExpression}.

\begin{definition}\label{Def-Fa}
Let $f$ be the polynomial whose zeta function we wish to compute,
and write
\[ X_{0}f  = \sum_{j \in J}\bar{a}_{j}X^{j}, \]
where $J$ is the support of $X_{0}f$. Let $a_{j}$ be the
Teichm\"{u}ller lifting of $\bar{a}_{j}$.
Let $F$ be the formal power series in the indeterminates $X$ with
coefficients in $\Rp$ given by
\[
F(X) = \prod_{j \in J} \theta(a_{j}X^{j}).
\]
Let $F^{(\extdeg)}(X)$ be the formal power series in the indeterminates
$X$ with coefficients in $\Rp$ given by
\[ F^{(\extdeg)}(X) =  \prod_{j \in J}
\prod_{s = 0}^{\extdeg -1} \theta((a_{j}X^{j})^{p^{s}}).\]
\end{definition}

The relation between $F(X)$ and $F^{(\extdeg)}(X)$ is clear.

\begin{lemma}\label{Lem-RelFFa}
Let the power series $F^{(\extdeg)}$ and $F$ be as in
Definition \ref{Def-Fa}. Then
\[ F^{(\extdeg)}(X) = \prod_{i = 0}^{\extdeg -1}\tau^{i}(F(X^{p^{i}}))
\]
where the map $\tau$ acts coefficient-wise on the
power series $F$.
\end{lemma}

The power series $F^{(\extdeg)}(X)$ relates to rational point counting
in the following way.

\begin{proposition}\label{Prop-AnalyticExpression}
Let $f \in \fq[X_{1},\dots,X_{n}]$ and $F^{(\extdeg)} \in
\Rp[[X]]$ be as in Definition \ref{Def-Fa}. Denoting by
$N_{k}^{*}$ the number of solutions to the equation $f = 0$ in the
affine torus $(\fqk^{*})^{n}$, we have
\[
q^{k}N_{k}^{*} - (q^{k} - 1)^{n} = \sum_{x^{q^{k}-1} = 1}
F^{(\extdeg)}(x)F^{(\extdeg)}(x^{q})\dots F^{(\extdeg)}(x^{q^{k-1}})\]
where the summation is over the Teichm\"{u}ller lifting
of points on the torus $(\fqk^{*})^{n+1}$.
\end{proposition}

\begin{proof}
For any point $\bar{x}$ in $(\fqk)^{n+1}$ with Teichm\"{u}ller
lifting $x$ we have that
\[
\begin{array}{rcl}
\Psi_{k}(\bar{x}_{0}f(\bar{x}_{1},\dots,\bar{x}_{n}))
 & = & \Psi_{k}(\sum_{j \in J} \bar{a}_{j}\bar{x}^{j})\\
 & = & \prod_{j \in J} \Psi_{k}(\bar{a}_{j}\bar{x}^{j})\\
 & = & \prod_{j \in J} \Phi_{k}(a_{j}x^{j})\\
            & = & \prod_{j \in J}\{\prod_{i = 0}^{\extdeg k-1}
            \theta((a_{j}x^{j})^{p^{i}})\}\\
            & = & \prod_{j \in J}\{\prod_{i = 0}^{k-1}
\prod_{s = 0}^{\extdeg -1} \theta((a_{j}x^{j})^{q^{i}p^{s}})\}\\
& = & \prod_{i = 0}^{k-1} \{\prod_{j \in J}
(\prod_{s = 0}^{\extdeg -1} \theta((a_{j}^{q^{i}}x^{jq^{i}})^{p^{s}}))\}\\
& = & \prod_{i = 0}^{k-1} \{\prod_{j \in J}
(\prod_{s = 0}^{\extdeg -1} \theta((a_{j}(x^{q^{i}})^{j})^{p^{s}}))\}\\
& = & F^{(\extdeg)}(x)F^{(\extdeg)}(x^{q})\dots F^{(\extdeg)}(x^{q^{k-1}}),
\end{array}
\]
where $F^{(\extdeg)}(X)$ is given by
\[ F^{(\extdeg)}(X) =  \prod_{j \in J}
\prod_{s = 0}^{\extdeg-1} \theta((a_{j}X^{j})^{p^{s}}).\]
(We pause to justify the steps above: the first four equalities
follow straight from definitions and
from the homomorphic property of $\Psi_{k}$; the fifth and sixth
by rearrangement; and the
seventh since $a_{j}$ satisfies $a_{j}^{q} = a_{j}$.)

Thus we have
\[ \sum_{\bar{x} \in (\sfqk^{*})^{n+1}}
\Psi_{k}(\bar{x}_{0}f(\bar{x}_{1},\dots,\bar{x}_{n})) =
\sum_{x^{q^{k} - 1} = 1} F^{(\extdeg)}(x)F^{(\extdeg)}(x^{q})
\dots F^{(\extdeg)}(x^{q^{k-1}})\] where the latter sum is over
the Teichm\"{u}ller lifting in $\Omega^{n+1}$ of points in
$(\fqk^{*})^{n+1}$. Combining this with Lemma \ref{Lem-ExpSum}
gives us the result.
\end{proof}

\subsection{Decay rates and weight functions}

We now describe the decay rates of the coefficients of the
power series $F^{(\extdeg)}$ and $F$.
Specifically, we obtain lower bounds for the
$p$-adic order of the coefficients of the power
series $F$ expressed in terms of a certain weight function
on integer vectors.

Write $F = \sum_{r} F_{r}X^{r}$ where the sum is over
non-negative integer vectors in $\Z_{\geq 0}^{n+1}$.
Let $A$ be the $(n+1) \times |J|$ matrix
whose columns are $j = (j_{0},j_{1},\dots,j_{n}) \in J$. Then
from Definition \ref{Def-Fa} one sees
\begin{equation}\label{Eqn-Frcoeffs}
F_{r} = \sum_{u}(\prod_{j \in J} \lambda_{u_{j}} a_{j}^{u_{j}})
\end{equation}
where the outer sum is over all $|J|$-tuples $u = (u_{j})$ of
non-negative integers such that
\begin{eqnarray}\label{Eqn-Au=r}
Au = r,
\end{eqnarray}
thinking of $u$ and $r$ as column vectors.
Since $j_0=1$ for all $(j_0,\dots, j_n)\in J$,
the first row of the matrix $A$ is the vector $(1,1,\dots, 1)$.
The first equation in the above linear system is then
\begin{eqnarray}\label{Eqn-r0}
\sum_{j\in J} u_j =r_0.
\end{eqnarray}
Now $F_{r}$ is
zero if (\ref{Eqn-Au=r}) has no solutions. Otherwise, since
$\ord(\lambda_{u_{j}} a_{j}^{u_{j}}) = \ord(\lambda_{u_{j}})$ from
(\ref{Eqn-DecayLambdas}),(\ref{Eqn-Frcoeffs}) and (\ref{Eqn-r0})
we get
\begin{eqnarray}\label{Eqn-Frdecay}
\ord(F_{r}) \geq \inf_{u} \left\{\sum_{j \in J}
\frac{(p-1)u_{j}}{p^{2}}
\right\} =\frac{p-1}{p^2} r_0,
\end{eqnarray}
where the inf is over all non-negative integer
vector solutions $u$ of (\ref{Eqn-Au=r}).
We now define a weight function $w$ such that
(\ref{Eqn-Frdecay}) gives
estimates on $\ord(F_{r})$ in terms of this weight function.

Let $\delta_{1} \subset \R^{n}$ denote the convex hull of the support
of $f$ (the set of exponents of non-zero terms). Let $\delta_{2}
\subset \R^{n}$ be the convex hull of the origin and the $n$ points
$(d,0,\dots,0),(0,d,\dots,0),
\dots,(0,\dots,0,d)$ where $d$ is the total degree of
$f$. We call $\delta_{1}$ the Newton
polytope of $f$; the polytope $\delta_{2}$ is just a simplex
containing $\delta_{1}$.

\begin{definition}\label{Def-Deltas}
Let $\delta$ be any convex polytope with
integer vertices such that $\delta_{1} \subseteq \delta \subseteq
\delta_{2}$. Denote by $\Newt$ the convex polytope in $\R^{n+1}$
obtained by embedding $\delta$ in $\R^{n+1}$ via the map
$x \mapsto (1,x)$ for $x \in \R^{n}$, and taking the convex hull with
the origin. Denote by $C(\Newt)$ the cone generated in $\R^{n+1}$
as the positive hull of $\Newt$.
So $C(\Newt)$ is the union of all rays emanating
from the origin and passing through $\Newt$.
\end{definition}

Ultimately, in Sections \ref{Sec-BitCompAux} and
\ref{Sec-FinalProof}, we shall only be interested in the simplest
choice of polytope $\delta = \delta_{2}$ (although the choice $\delta =
\delta_{1}$ leads to the most refined algorithm).
Letting
$\Newt_{1}$ denote the polytope in $\R^{n+1}$ obtained by
choosing $\delta = \delta_{1}$ we see that $C(\Newt_{1})$ is the
cone generated by the exponents of non-zero terms in $X_{0}f$.
Equation (\ref{Eqn-Au=r}) has no non-negative
integer (or even real) solutions when
$r$ does not lie in $C(\Newt_{1})$. Thus for any
choice of $\delta (\supseteq \delta_{1})$ and corresponding
$\Newt (\supseteq \Newt_{1})$,
all exponents of $F$ lie in the cone $C(\Newt)$.

\begin{definition}\label{Def-wf}
Define a weight function $w$ from $\R^{n+1}$ to
$\R \cup \{\infty\}$ in the following
way: For $r = (r_{0},r_{1},\dots,r_{n}) \in \R^{n+1}$ define
\[
w(r) = \left\{
\begin{array}{rl}
r_{0} & \mbox{if } r\in C(\Newt),\\
\infty & \mbox{otherwise.}
\end{array}
\right.
\]
In particular
$w(r)$ is a non-negative integer
for any $r \in C(\Newt)\cap \Z^{n+1}$.
\end{definition}

\begin{note}
Choosing $\delta = \delta_{1}$
corresponds to working with
the weight function of Adolphson and Sperber \cite{AS87}, and taking
$\delta = \delta_{2}$ to Dwork's original weight function \cite{BD60}.
Note that Dwork's weight function can equivalently
be defined as $w(r) = r_{0}$ if $r_{1} + \dots + r_{n} \leq r_{0}d$ and
$\infty$ otherwise, where $d$ is the total degree of $f$.
\end{note}

The weight function has a simple
geometric interpretation: For a real number
$c$ define $c\Newt = \{cx\,|\,x \in \Newt\}$.
The next lemma is straightforward.

\begin{lemma}\label{Lem-GeomPropwr}
For any point $r \in C(\Newt)$ we have that $w(r)$ is the smallest
non-negative number $c$ such that $r \in c\Newt$. If
$r \not \in C(\Newt)$ then $w(r) = \infty$.
\end{lemma}

By (\ref{Eqn-Frdecay}) and the sentence preceding
Definition \ref{Def-wf}, we have the following result.

\begin{lemma}\label{Lem-Frinequalities}
We have the inequality
\[\ord(F_{r}) \geq w(r)\frac{p-1}{p^{2}}.\]
\end{lemma}

We also note the following simple property whose
proof is straightforward.

\begin{lemma}\label{Lem-Propwqwb}
Let $r, r^{\prime} \in \R^{n+1}$ and $k$ a non-negative integer.
Then $w(kr) = kw(r)$ and $w(r + r^{\prime}) \leq
w(r) + w(r^{\prime})$. In particular when $w(r^{\prime}) \ne \infty$,
$$w(kr-r^{\prime})\geq kw(r)-w(r^{\prime}).$$
\end{lemma}

We shall work in certain subrings of $\Rp[[X]]$ defined
in terms of the weight function.

\begin{definition}\label{Def-Lw}
Define $L_{\Delta}$ to be the subring of
$\Rp[[X]]$ given by
\[ L_{\Delta} = \{\sum_{r\in C(\Newt)\cap \sZ^{n+1}} A_{r}X^{r}\,|\,A_{r}\in R \}.\]
\end{definition}

Thus, $L_{\Delta}$ is just the ring of all power series over $R$ whose
terms have exponents in the cone $C(\Newt)$.
Certainly $F \in L_{\Delta}$ and from Lemma \ref{Lem-RelFFa} we see
easily that $F^{(\extdeg)} \in L_{\Delta}$.

\begin{lemma}\label{Lem-FFainLw}
The power series $F$ and $F^{(\extdeg)}$ both belong to $L_{\Delta}$.
\end{lemma}

This concludes all results in this section which shall be essential
to the proof of our modular version of Dwork's Trace Formula.
We conclude with a definition and some comments which
we will refer to in the analysis of the running time of our algorithm.

\begin{definition}\label{Def-Lb}
For any positive real number $b$,
we define a set of power series by
\[ L_{\Delta}(b) = \{\sum_{r} A_{r}X^{r} \in L_{\Delta}\,|\,
\ord{A_{r}} \geq b w(r)\}.\]
\end{definition}

The set $L_{\Delta}(b)$ is
easily seen to be a subring of $L_{\Delta}$.
Elements in $L_{\Delta}(b)$ for large $b$ can be thought of as having
fast decaying coefficients.
Such rings will reduce to rings of small
dimension modulo small powers of $p$.

We have that
\begin{eqnarray}
F & \in & L_{\Delta}\left(\frac{p-1}{p^{2}}\right)\label{Eqn-FLw}\\
F^{(\extdeg)} & \in & L_{\Delta}\left(\frac{p-1}{qp}\right).
\end{eqnarray}
The first is immediate from the inequality
in Lemma \ref{Lem-Frinequalities}
and the second follows
since $\tau^{i}(F(X^{p^{i}})) \in L_{\Delta}(\frac{p-1}{p^{i}p^{2}})$ for each
$0 \leq i \leq \extdeg -1$. Thus for $\extdeg > 1$
the coefficients of $F^{(\extdeg)}$ decay more slowly than those of
$F$ itself.

\section{Dwork's trace formula} \label{Sec-DTF}

\subsection{Lifting Frobenius}

We now introduce Dwork's ``left inverse of
Frobenius'' mapping $\psi_{p}$ on the ring $\Rp[[X]]$.

\begin{definition}\label{Def-Psi}
Let $\psi_{p}$ be defined on the monomials in $\Rp[[X]]$ by
\[
\psi_{p}(X^{r}) = \left\{
\begin{array}{ll}
X^{r/p} & \mbox{if }p|r\\
0 & \mbox{otherwise}
\end{array}
\right.
\]
and extend $\psi_{p}$ by $\tau^{-1}$-linearity to all of
$\Rp[[X]]$. That is
\[ \psi_{p}(\sum_{r} A_{r} X^{r}) = \sum_{r} \tau^{-1}(A_{r})
\psi_{p}(X^{r}) = \sum_{r} \tau^{-1}(A_{pr}) X^{r}.\]
\end{definition}

Here $p|r$ means that $p$ divides all of the entries in the integer
vector $r$. This map is a left inverse of the ``Frobenius''
map on the ring $\Rp[[X]]$ which takes a power series
$\sum_{r}A_{r}X^{r}$ to $\sum_{r}\tau(A_{r})X^{pr}$.

\begin{definition}\label{Def-alphas}
Let $\alpha_{\extdeg}$ be the map from $\Rp[[X]]$ to itself defined as
\[ \alpha_{\extdeg} = \psi_{p}^{\extdeg} \circ F^{(\extdeg)}.\]
Precisely, this is the
map which is the composition of multiplication by the power
series $F^{(\extdeg)}$ followed
by the mapping $\psi_{p}^{\extdeg}$ on the ring $\Rp[[X]]$.
(Notice that $\psi_{p}^{\extdeg}$
just acts as
\[ \psi_{p}^{\extdeg}(\sum_{r} A_{r} X^{r}) = \sum_{r} A_{qr} X^{r}\]
since $\tau^{-\extdeg}(A_{r}) = A_{r}$.) Let the map $\alpha$ from
$\Rp[[X]]$ to itself be defined as
\[\alpha = \psi_{p} \circ F.\]
Thus $\alpha$ is multiplication by $F$ followed by the mapping
$\psi_{p}$.
\end{definition}

We have the following result relating these two maps,
which shall be of crucial importance in our derivation of an efficient
algorithm.

\begin{lemma}\label{Lem-alphapower}
With $\alpha_{\extdeg}$ and $\alpha$ as in Definition \ref{Def-alphas} we have
\[\alpha_{\extdeg} = \alpha^{\extdeg}\]
where the second exponent is a power
under composition.
\end{lemma}

\begin{proof}
Firstly let $H \in \Rp[[X]]$ and denote by $\psi_{p} \circ
H(X^{p})$ the map composed of multiplication by $H(X^{p})$ followed
by $\psi_{p}$. We claim that
\begin{equation}\label{Eqn-Proppsi}
\psi_{p} \circ H(X^{p}) = \tau^{-1}(H(X))
 \circ \psi_{p}
\end{equation}
To see this write $H = \sum_{r}H_{r}X^{pr}$. Then we have that
\[\psi_{p} \circ H(X^{p}) = \sum_{r} \tau^{-1}(H_{r})(\psi_{p} \circ
X^{pr}).\]
Here the infinite series is interpreted as a mapping. Now
$\psi_{p} \circ X^{pr} = X^{r} \circ \psi_{p}$ as these two maps are
$\tau^{-1}$-linear and agree on monomials. Hence we have $\psi_{p} \circ
H(X^{p}) = \sum_{r}\tau^{-1}(H_{r})(X^{r} \circ \psi_{p}) =
\tau^{-1}(H(X)) \circ \psi_{p}$.

Next we claim that for any $b \geq 1$ and power series
$H$ we have
\[ \psi_{p}^{b} \circ \prod_{i = 0}^{b-1} \tau^{i}(H(X^{p^{i}}))
 = (\psi_{p} \circ H(X))^{b}.\]
We prove this by induction, the result trivially holding if
$b = 1$. For $b > 1$, by $b-1$ applications of (\ref{Eqn-Proppsi}) we get
\[ \psi_{p}^{b} \circ \prod_{i = 0}^{b-1} \tau^{i}(H(X^{p^{i}}))
= (\psi_{p} \circ H(X)) \circ (\psi_{p}^{b-1} \circ \prod_{i = 0}
^{b-2} \tau^{i}(H(X^{p^{i}})).\]
The second claim then follows by induction. Putting $H = F$ and
$b = \extdeg$ we get the required result.
\end{proof}

The map $\alpha_{\extdeg}$ is linear and continuous, in the sense that
\[\alpha_{\extdeg}
(\sum_{r}A_{r}X^{r}) = \sum_{r} A_{r} \alpha_{\extdeg}(X^{r})\]
for any element $\sum_{r} A_{r}X^{r} \in \Rp[[X]]$.
The map $\alpha$ is $\tau^{-1}$-linear and continuous, in
the sense that
\[\alpha(\sum_{r} A_{r}X^{r}) = \sum_{r}
\tau^{-1}(A_{r})\alpha (X^{r}).\]

From Lemma \ref{Lem-FFainLw} the next lemma follows easily.

\begin{lemma}\label{Lem-Lwstable}
The subring $L_{\Delta}$ is stable under both
$\alpha$ and $\alpha_{\extdeg}$.
\end{lemma}

Both maps when restricted to the subring $L_{\Delta}$ are determined
by their action on the monomials $X^{r}$ (with $w(r) < \infty$),
which we now consider.

\subsection{Matrix representations of mappings}

Recall that $C(\Newt)$ is the cone in $\R^{n+1}$ from Definition
\ref{Def-Deltas}.
The set
\[
\Gamma_{\Delta}= \{X^{u}\,|\,u \in C(\Newt) \cap \Z_{\geq 0}^{n+1} \}
\]
written as a row vector, is a {\it formal basis} for the space $L_{\Delta}$. Precisely, this
means that any power series in $L_{\Delta}$ may be written in exactly
one way as an infinite sum $\sum_{u} A_{u}X^{u}$ with $X^{u}$ in the
above set. Notice that this is different from the usual notion of a
basis in linear algebra, since we allow infinite combinations of basis
elements. It is also different from the notion of an ``orthonormal
basis'' in the literature, where one requires that the coefficient
$A_u$ goes to zero as $w(u)$ goes to $\infty$ (see \cite{DW1}
for a more detailed discussion of these notions).

By Lemma \ref{Lem-Lwstable}, both $\alpha$ and $\alpha_{\extdeg}$
send the ring $L_{\Delta}$ to itself. We define certain matrices
associated to the maps $\alpha$ and $\alpha_{\extdeg}$ restricted to
$L_{\Delta}$ with regard to the formal row basis $\Gamma_{\Delta}$ of monomials.

\begin{definition}\label{Def-MMa}
Let the infinite matrices $M$ and $M_{\extdeg}$
have columns describing the images of the monomials
$X^{v} \in L_{\Delta}$ under the maps $\alpha$ and $\alpha_{\extdeg}$
with respect to our formal row basis $\Gamma_{\Delta}$:
\[ \alpha(\Gamma_{\Delta})=\Gamma_{\Delta}M, \
\alpha_{\extdeg}(\Gamma_{\Delta})=\Gamma_{\Delta}M_{\extdeg}.
\]
Specifically, the $(u,v)$th entries of
$M$ and $M_{\extdeg}$ for $u,v \in C(\Newt)$
are $m_{uv}$  and $m^{(\extdeg)}_{uv}$ respectively where
\[
\begin{array}{rcl}
m_{uv} & = & \tau^{-1}(F_{pu-v})\\
m^{(\extdeg)}_{uv} & = & F^{(\extdeg)}_{qu-v}.
\end{array}
\]
Here $F^{(\extdeg)} = \sum_{r} F^{(\extdeg)}_{r}X^{r}$ and as before
$F = \sum_{r} F_{r}X^{r}$, and we take the coefficients of exponents
$r$ with negative entries to be zero.
\end{definition}

(Note that as yet we have not ordered the basis; however,
we shall choose a convenient ordering in the proof of Theorem
\ref{Thm-MDTF}.)
We have by Lemmas \ref{Lem-Frinequalities} and \ref{Lem-Propwqwb}
\begin{eqnarray}
    \ord(F_{pu-v}) & \geq & \frac{p-1}{p^{2}}w(pu-v)\nonumber\\
    & \geq & \frac{p-1}{p}\left(w(u) - \frac{1}{p}w(v)\right)\label{Eqn-Mdecay}
\end{eqnarray}
and certainly $\ord(m_{uv}) = \ord(\tau^{-1}(F_{pu-v})) =
\ord(F_{pu-v})$.

Notice that the matrix powers $M_{\extdeg}^{k}$ and $M^{k}$
are defined for every positive integer $k$
since the entries in $M_{\extdeg}^{k}$, say, are just finite sums of
the entries in $M_{\extdeg}$ and $M_{\extdeg}^{k-1}$.
This follows for $M_{\extdeg}$, say,
since all entries $m_{uv}^{(\extdeg)}$ in $M_{\extdeg}$ are
zero when the vector $qu-v$ contain negative
entries.
We define the trace of an infinite matrix to be the sum of its
diagonal entries, when this sum converges, and $\infty$ when
the sum does not converge.
We shall see shortly that the trace of the infinite matrix $M_{\extdeg}^{k}$
is finite. We write this as $\Tr(M_{\extdeg}^{k})$.

We can now prove the following result of Dwork.

\begin{theorem}[Dwork's Trace Formula]\label{Thm-DTF}
Let $f \in \fq[X_{1},\dots,X_{n}]$ and let $\alpha_{\extdeg}$
be the mapping on the ring $\Rp[[X]]$ given as
$\alpha_{\extdeg} = \psi_{p}^{\extdeg} \circ F^{(\extdeg)}$, where
$F^{(\extdeg)}$ and $\psi_{p}$ are described in Definitions
\ref{Def-Fa} and \ref{Def-Psi}.
Let $M_{\extdeg}$ denote the infinite matrix representing the
map $\alpha_{\extdeg}$ restricted to the subring
$L_{\Delta}$ as described in Definition
\ref{Def-MMa}.
For $k \geq 1$,
denote by $N_{k}^{*}$ the number of solutions to the
equation $f = 0$ in the torus $(\fqk^{*})^{n}$. Then
\[
(q^{k} - 1)^{n+1} \Tr(M_{\extdeg}^{k}) = q^{k}N^{*}_{k} - (q^{k} - 1)^{n}.
\]
\end{theorem}

\begin{proof}
(The following is a hybrid of the
matrix proof given in \cite{DW96} and
Dwork's original argument \cite{BD60}.)

By Proposition \ref{Prop-AnalyticExpression} we have that
\[
\begin{array}{rcl}
q^{k}N^{*}_{k} - (q^{k} - 1)^{n}
& = & \sum_{x^{q^{k} - 1} = 1} F^{(\extdeg)}(x)F^{(\extdeg)}(x^{q})\dots
F^{(\extdeg)}(x^{q^{k-1}})
\end{array}
\]
where the sum is over all $(n+1)$-tuples of $(q^{k} - 1)$st roots
of unity in $\Omega$ (namely the Teichm\"{u}ller lifting in
$\Omega^{n+1}$ of points on the torus $(\fqk^{*})^{n+1}$).

We first consider the case $k = 1$. Since $F^{(\extdeg)} \in L_{\Delta}$,
we can write
$F^{(\extdeg)}(X) = \sum_{r} F^{(\extdeg)}_{r}X^{r}$
where the sum is over all lattice vectors $r$
which belong to $C(\Newt)$.
Then the latter sum is
\[
\begin{array}{rcl}
\sum_{x^{q - 1} = 1} F^{(\extdeg)}(x) & = & \sum_{x^{q-1} = 1} (\sum_{r}
F^{(\extdeg)}_{r}x^{r})\\
& = & \sum_{r}F^{(\extdeg)}_{r}(\sum_{x^{q-1} = 1} x^{r})\\
& = & (q-1)^{n+1}\sum_{r, (q-1)|r} F^{(\extdeg)}_{r}\\
& = & (q-1)^{n+1}\sum_{s} F^{(\extdeg)}_{(q-1)s}\\
& = & (q-1)^{n+1}\Tr(M_{\extdeg}).
\end{array}
\]
Here by $(q-1)|r$ we mean that the integer $q-1$ divides every
entry in the vector $r$. Also, we use the fact that for any
integer $r_{i}$ \cite[page 120]{NK}
\[ \sum_{x_{i}^{q-1} = 1} x_{i}^{r_{i}} =
\left\{
\begin{array}{ll}
q - 1 & \mbox{if }(q-1)|r_{i}\\
0 & \mbox{otherwise.}
\end{array}
\right.
\]

For $k > 1$, define $M_{\extdeg k}$ to be the matrix for the map
$\psi_{p}^{\extdeg k} \circ \prod_{i = 0}^{k-1}F^{(\extdeg)}
(X^{q^{i}})$ with respect to
the formal basis $\Gamma_{\Delta}$. Then by an analogous argument to
that in the case $k = 1$ we see
\[ \sum_{x^{q^{k} - 1} = 1} \prod_{i = 0}^{k-1} F^{(\extdeg)}(x^{q^{i}}) =
(q^{k} - 1)^{n+1}\Tr(M_{\extdeg k}).\]
Using (\ref{Eqn-Proppsi}) in the proof of Lemma \ref{Lem-alphapower}
one sees that
\[ \psi_{p}^{\extdeg k} \circ \prod_{i =
0}^{k-1}F^{(\extdeg)}(X^{q^{k}}) = (\psi_{p}^{\extdeg} \circ
F^{(a)}(X))^{k}.\]
Since both maps are linear it follows that $M_{\extdeg k} =
M_{\extdeg}^{k}$, and the theorem is proved.
\end{proof}

To compute the trace of the matrix in Dwork's formula we shall use
the following matrix identity derived from Lemma \ref{Lem-alphapower}.

\begin{lemma}\label{Lem-Matrices}
Let $M_{\extdeg}$ and $M$ denote the matrices for the maps
$\alpha_{\extdeg}$ and
$\alpha$ described above. Then
\[
M_{\extdeg} = \prod_{i = 0}^{\extdeg-1} \tau^{-i}(M)
= M M^{\tau^{-1}}\cdots M^{\tau^{-(\extdeg-1)}}
\]
where the map $\tau^{-i}$ acts entry wise on the matrix $M$.
\end{lemma}

\begin{proof}
First suppose that $N_{1}$ and $N_{2}$ are matrices representing
maps $\beta_{1}$ and $\beta_{2}$ on some subspace of $\Rp[[X]]$.
We assume that $\beta_{1}$ is
$\tau^{-1}$-linear and $\beta_{2}$ is $\tau^{-j}$-linear for some
$j \in \Z$. Then it is not
difficult to prove that the matrix for the $\tau^{-(j+1)}$-linear
map $\beta_{1} \circ \beta_{2}$ is just $N_{1}\tau^{-1}(N_{2})$.

By Lemma \ref{Lem-alphapower} we have
$\alpha_{\extdeg} = \alpha^{\extdeg}$. We claim that for any positive
integer $b$ the matrix for the map $\alpha^{b}$ is
$\prod_{i = 0}^{b-1}\tau^{-i}(M)$. The result is
trivially true if $b = 1$. For $b > 1$ we have
$\alpha^{b} = \alpha \circ \alpha^{b-1}$. By induction
the matrix for $\alpha^{b-1}$ is $\prod_{i = 0}^{b-2}\tau^{-i}(M)$.
By $\tau^{-1}$-linearity of $\alpha$ it follows from the
observations in the preceding paragraph that
the matrix for $\alpha^{b}$ is $M \tau^{-1}(\prod_{i = 0}^{b-2}
\tau^{-i}(M)) = \prod_{i = 0}^{b-1} \tau^{-i}(M)$. The required
result now follows by taking $b = \extdeg$.
\end{proof}

We note in passing that $M_{\extdeg}$ in Lemma \ref{Lem-Matrices}
also equals
\[
\tau^{\extdeg -1}(S)\tau^{\extdeg -2}(S)\dots S,
\]
where $
S = \tau(M)$ is the matrix with $(u,v)$th entry simply
$F_{pu-v}$. Although this is slightly more desirable from
a practical point of view we shall not use this expression
for $M_{\extdeg}$.

\subsection{Modular reduction of the trace formula}\label{Sec-MDTF}

We now examine the reduction of the Dwork trace formula modulo
a power $p^{N}$ of $p$. Observe that both sides in the
trace formula, and all the entries in the matrices $M$ and
$M_{\extdeg}$, are elements of $\Rp$, and so this reduction is defined.

We first recall some notation: Let $C(\Newt)$ be the cone in
$\R^{n+1}$ from Definition \ref{Def-Deltas}, and
$L_{\Delta}$ denote the ring of power series over $\Rp$ whose monomials
have exponents lying in $C(\Newt)$ (Definition \ref{Def-Lw}).
Let $M$ denote the matrix for the $\tau^{-1}$-linear map
$\alpha = \psi_{p} \circ F$ with respect to the formal row basis $\Gamma_{\Delta}$
(Definition \ref{Def-MMa}). Here
$\psi_{p}$ is the ``left inverse of Frobenius'' given
in Definition \ref{Def-Psi}
and $F$ is the power series obtained from the
polynomial $f$ as in Definition \ref{Def-Fa}.

\begin{theorem}\label{Thm-MDTF}
Let $N$ denote any positive
integer and $A_N$ the finite square matrix over the finite ring
$\Rp/(p^{N})$ obtained by reducing modulo $p^{N}$
all those entries in $M$ whose
rows and columns are indexed by vectors
$u \in C(\Newt) \cap \Z_{\geq 0}^{n+1}$ with
$w(u) < (p/(p-1))^{2}N$.
Then
\[ (q^{k} - 1)^{n+1}
\Tr((\prod_{i = 0}^{\extdeg -1}\tau^{-i}(A_N))^{k}) = q^{k}N^{*}_{k} - (q^{k} -
1)^{n} \bmod{p^{N}}\]
where $N^{*}_{k}$ is the number of solutions to the equation
$f = 0$ in the affine torus $(\fqk^{*})^{n}$.
Moreover, the size of
$A_{N}$ is $W = \# (t\Newt)$, where $t = \lceil
p^{2}N/(p-1)^{2} \rceil - 1$ and $\#(t\Newt)$ is
number of lattice points in a dilation by a factor $t$ of the
polytope $\Newt$.
\end{theorem}

\begin{proof}
For any finite or infinite matrix $L$
with coefficients
in $\Rp$,
we define
$\overline{L}$ to be the matrix obtained by reducing
all its entries modulo
$p^{N}$. Thus $\overline{L}$ has entries in
$\Rp/(p^{N})$.

The theorem will follow from the
Dwork Trace Formula (Theorem \ref{Thm-DTF}) once we find
a suitable expression for the reduction modulo $p^{N}$ of the
trace of the matrix $M_{\extdeg}^{k}$. This is equal to
the trace of the matrix $\overline{M_{\extdeg}^{k}} =
(\overline{M_{\extdeg}})^{k}$.
By Lemma \ref{Lem-Matrices} this matrix
can be computed as a matrix product from the matrix
$\overline{M}$.

By inequality
(\ref{Eqn-Mdecay})
every entry $m_{uv} = \tau^{-1}(F_{pu-v})$ in $M$ satisfies
\[
\ord(m_{uv}) \geq \frac{p-1}{p}\left(w(u) - \frac{w(v)}{p} \right).
\]
Define $t$ to be the greatest integer less than $(p/(p-1))^{2}N$.
When $w(u) \geq w(v)$ and $w(u) > t$ we have
\begin{equation}\label{Eqn-MatrixEntries}
\frac{p-1}{p}\left(w(u) - \frac{w(v)}{p} \right)
\geq \frac{(p-1)^2}{p^2}w(u)
\geq \frac{(p-1)^2}{p^2}(t+1) \geq N.
\end{equation}
Recall that we have not yet ordered the basis. Now
choose any total ordering on the basis set such that for
distinct lattice points $u,v \in C(\Newt)$, the monomial $X^{v}$
comes before $X^{u}$ if $w(v) < w(u)$.
Thus,
by the inequalities in
(\ref{Eqn-MatrixEntries}) and the choice of ordering of the basis,
the matrix $\overline{M}$
is of the form
\begin{equation}\label{Eqn-ReducedMatrix}
\left(
\begin{array}{cc}
    A_N  & B_N\\
    0  & C_N
\end{array}
\right)
\end{equation}
where $C_N$ is a strictly upper triangular infinite matrix and
$A_N$ is the finite square matrix indexed by lattice points
$u \in C(\Newt)$ such that $w(u) \leq t$.
The size of $A_N$ is the number of lattice points $u$ with
$w(u) \leq t$.
By Lemma \ref{Lem-GeomPropwr}
this is exactly the number of lattice points in the
polytope $t\Newt$.

By Lemma \ref{Lem-Matrices} and modular reduction,  one has
\[ (\overline{M_{\extdeg}})^{k} = (\prod_{i = 0}^{\extdeg-1}
\tau^{-i}(\overline{M}))^{k}.\]
By (\ref{Eqn-ReducedMatrix})
we see that $(\overline{M_{\extdeg}})^{k}$ is of the form
\[
\left(
\begin{array}{cc}
(\prod_{i = 0}^{\extdeg -1} \tau^{-i}(A_N))^{k} & B_N^{\prime}\\
0 & C_N^{\prime}
\end{array}
\right)
\]
where $C_N^{\prime} = (\prod_{i = 0}^{\extdeg -1} \tau^{-i}(C_N))^{k}$
is strictly upper triangular.

Hence the trace of $(\overline{M_{\extdeg}})^{k}$ equals the trace
of the finite matrix
\begin{equation}\label{Eqn-MatrixProduct}
(\prod_{i = 0}^{\extdeg -1}\tau^{-i}(A_N))^{k}.
\end{equation}

The theorem now follows from (\ref{Eqn-MatrixProduct}) and Theorem
\ref{Thm-DTF}.
\end{proof}

\section{Algorithms}\label{Sec-ALG}

In this section we present an algorithm for counting points based
upon Theorem \ref{Thm-MDTF}, and complete the proofs of the
results in the introduction. This is a relatively straightforward
matter, although the precise complexity estimates require
a little care.

\subsection{Toric point counting algorithm}

We first give the algorithm.

\begin{algorithm}\label{Alg-Torus}
{\bf Toric point counting}\\
Input: Positive integers $\extdeg,k,n,d$ and a prime $p$; a polynomial $f$; a
polytope $\Newt$.
(Here $f$ is a
polynomial in $n$ variables of total degree $d$ with
coefficients in the field $\fq$ where $q = p^{\extdeg}$. The polytope $\Newt$
is as in Definition \ref{Def-Deltas}. We assume a
model of $\fq$ is given as in Section \ref{Sec-ContLiftFrob}.)\\
Output: The number of solutions $N_{k}^{*}$ to the equation $f=0$ in
the torus $(\fqk^{*})^{n}$.\\

\hspace{-\parindent}Step 0: Set $N = (n+1)\extdeg k$, where $q =
p^{\extdeg}$.\\

\hspace{-\parindent}Step 1: Compute the polynomial $F \bmod{p^{N}}$
in the ring $(\Rp/(p^{N}))[X]$ where $F$ is the power series
in Definition \ref{Def-Fa}.
Specifically, writing $X_{0}f = \sum_{j \in
J} \bar{a}_{j}X^{j}$ we have $F = \prod_{j \in J}
\theta (a_{j}X^{j}) \bmod{p^{N}}$. Here $a_{j}$ is the Teichm\"{u}ller
lifting of $\bar{a}_{j}$ and $\theta(z) = \exp(\pi(z -
z^{p}))$. (See Section \ref{Sec-PadicTheory} for a description of the
ring $\Rp/(p^{N})$ and the element $\pi$.)
\\

\hspace{-\parindent}Step 2: Construct the matrix $A_N$ which occurs in
the statement of Theorem \ref{Thm-MDTF}. Specifically, the matrix
$A_{N}$ is indexed by pairs $(u,v)$ where $u$ and $v$ are lattice points
in the dilation by a factor $t$ of the polytope $\Newt$, and
$t = \lceil p^{2} N/(p-1)^{2} \rceil - 1$. The $(u,v)$th
entry of $A_{N}$ is $\tau^{-1}$ of the coefficient
of $X^{pu-v}$ in the polynomial $F \bmod{p^{N}}$.
The action of $\tau$ is as described in Section \ref{Sec-PadicTheory}.
\\

\hspace{-\parindent}Step 3: Compute the product
$(\prod_{i = 0}^{\extdeg -1}\tau^{-i}(A_N))^{k}$.
Let $T$ denote the trace of this product.\\

\hspace{-\parindent}Step 4: Output
\[N_{k}^{*} = q^{-k}[((q^{k}-1)^{n+1}T + (q^{k} - 1)^{n}) {\rm mod}~p^N],\]
where the square brackets denote the smallest non-negative residue
modulo $p^N$.
\end{algorithm}

In the algorithm we assume that the polynomial is presented as input
explicitly via its list of non-zero terms. The manner of presentation
of $\Newt$ only affects the time required to find all lattice points
in the dilated polytope $t\Newt$. For concreteness, let us say it is
presented via its list of vertices, although any other reasonable
presentation would suffice.

\subsection{Proof of correctness of the algorithm}

We know by Theorem \ref{Thm-MDTF} that
\[ q^{k}N_{k}^{*} - (q^{k} - 1)^{n} =  (q^{k} - 1)^{n+1}\Tr((\prod_{i =
0}^{\extdeg -1}\tau^{-i}(A_N))^{k}) \bmod{p^{N}}\]
in the ring $\Rp/(p^{N})$.
Thus,
\[ q^{k}N_{k}^{*} = (q^{k} - 1)^{n}+ (q^{k} - 1)^{n+1}T \bmod{p^{N}}.\]
The lefthand side is a non-negative integer. The first term on the righthand side
is an integer. The second term on the righthand side is the reduction modulo $p^{N}$ of the
trace of $M_{\extdeg}^{k}$, which is known to be an integer by
Dwork's Trace Formula (Theorem \ref{Thm-DTF}), and thus it
is also an integer. Since
$N_{k}^{*} \leq (q^{k} - 1)^{n}$ it follows that the lefthand side is
smaller than $q^{k(n+1)}$. Hence in the case
$N \geq \extdeg k(n+1)$ we must have
\[q^kN_{k}^{*}= [((q^{k}-1)^{n+1}T + (q^{k} - 1)^{n}) {\rm mod}~p^N].\]
The proof is complete.

\subsection{Complexity analysis}\label{Sec-ComplexityAnalysis}

We shall use Big-Oh and Soft-Oh notation in our analysis of the
complexity of the above algorithm. If
$C_{1}$ and $C_{2}$ are real functions
we write $C_{1} = \Oh(C_{2})$ if
$|C_{1}| \leq c(|C_{2}|+1)$ for some positive constant $c$.
We write $C_{1} = \SoftOh(C_{2})$ if
$C_{1} = \Oh(C_{2}\log(|C_{2}| + 1)^{c^{\prime}})$ for some
constant $c^{\prime}$. Thus in the latter notation one ignores logarithmic
factors.

\subsubsection{Ring operations}\label{Sec-RO}

We shall first of all count the number of operations in the ring
$\Rp/(p^{N})$ required in Steps 1, 2, 3, ignoring for the time being
any other auxiliary computations. More precisely,
because of the complexity bounds in Lemma \ref{Lem-R0ops}
it is convenient for
our analysis to define a ``ring operation'' to be either
arithmetic or the evaluation of the map $\tau^{i}$, for
$1 \leq i \leq \extdeg - 1$, in the ring
$\Rp/(p^{N})$ (excluding precomputation).
In Section \ref{Sec-BitCompRing} we shall add back
in the small contribution from computing Teichm\"{u}ller liftings and
also the precomputation required for the maps $\tau^{i}$.
Similarly, in Section \ref{Sec-BitCompAux} we shall account for
the remaining operations in Steps 1, 2, 3,
arising mainly from computations with the
exponents of polynomials (at this stage we will restrict the
input polytope $\Newt$ to avoid complications from convex geometry).
The contributions from Steps 0 and 4 are easily seen to be absorbed
into the other estimates, and we shall not mention them again.

Our running time will be in terms of the
parameters $t,\tilde{t},W,\tilde{W}$. Here
\[
\begin{array}{ll}
 t = \lceil \left( \frac{p}{p-1} \right)^{2} N \rceil - 1,  & W =
\#(t\Newt),\\
 \tilde{t} = \lceil \frac{p^{2}}{p-1} N \rceil - 1,  & \tilde{W} =
\#(\tilde{t}\Newt),
\end{array}
\]
where $N = \extdeg k(n+1)$ and the $\#$ operator counts lattice
points in convex sets. The sizes of the sets $W$ and $\tilde{W}$
are the number of lattice points in certain
``truncated'' cones.
Since $\tilde{t}$ is
about $p$ times as large as $t$, the integer $\tilde{W}$ will be
around $p^{n+1}$ times as large as $W$, since we are working
in $n+1$ dimensional space. The integer $W$ is precisely the
size of the matrix which occurs in Step 2. The integer
$\tilde{W}$ will turn out to be the maximum number of
terms in the polynomial we compute in Step 1.
Define $L_{\Delta}((p-1)/p^{2}) \bmod{p^{N}}$ to be the ring
of polynomials obtained by reducing the coefficients of power series
in $L_{\Delta}((p-1)/p^{2}) \subseteq \Rp[[X]]$ modulo $p^{N}$.
Then
$\tilde{W}$ is the number of monomials which occur in the
finite ring $L_{\Delta}((p-1)/p^{2}) \bmod{p^{N}}$.

For Step 1 we have the following estimate.

\begin{lemma}\label{Lem-ComputingF}
Let $F$ be the power series given in Definition \ref{Def-Fa} and
$N$ a positive integer. The polynomial $F \bmod{p^{N}}$ may be
computed in
\[\Oh(|J|\tilde{W}^{2})\]
operations in the ring $\Rp/(p^{N})$.
\end{lemma}

\begin{proof}
By (\ref{Eqn-DecayLambdas}) we have that $\theta(z) \bmod{p^{N}}$ is
a polynomial of degree not greater than $p^{2}N/(p-1)$. Thus we can obtain
$\theta(z)$ via the formula
\[ \theta(z) = \exp(\pi z) \exp(-\pi z^{p})\]
by computing the first $\Oh(pN)$ terms in the expansion for $\exp(\pi z)$,
substituting $z = -z^p$, and one multiplication
of polynomials. Note that $\exp(\pi z)$ has $p$-adic integral
coefficients.
Thus $\theta(z)$ can be found in time
$\Oh((pN)^{2})$ operations in the ring $\Rp/(p^{N})$ using standard
polynomial arithmetic.

Now by (\ref{Eqn-FLw}) we have that $F \in L_{\Delta}((p-1)/p^{2})$.
Thus $F \bmod{p^{N}} \in L_{\Delta}((p-1)/p^{2}) \bmod{p^{N}}$.
One may then compute $F \bmod{p^{N}}$ directly
from Definition \ref{Def-Fa} in $|J|-1$ multiplications of polynomials of
the form $\theta(a_{j}X^{j}) \bmod{p^{N}}$.
Each such polynomial lies in the ring $L_{\Delta}((p-1)/p^{2})
\bmod{p^{N}}$, because
$\ord{(a_{j})} = 0, w(j) = 1$ and the coefficients of $\theta$ decay
at a suitable rate. Hence all computations required in computing
$F \bmod{p^{N}}$ involve polynomials in this ring.
Such polynomials have at most $\tilde{W}$ terms.
Exactly $|J|-1$ multiplications are required.
Thus the complexity is $\Oh(|J|\tilde{W^{2}})$ ring operations.
Noting that $pN = O(\tilde{W})$ we have the
result.
\end{proof}

\begin{note}
It is crucial here that we only need to compute $F \bmod{p^{N}}$ and
not $F^{(\extdeg)} \bmod{p^{N}}$, as one might attempt to do using a more
naive approach. The latter polynomial has very high degree ($\Oh(q)$)
because of the slow decay rate of the coefficients of $F^{(\extdeg)}$.
\end{note}

With regard to Step 2, given that the polynomial $F \bmod{p^{N}}$ has
already been computed the only task required is to identify those
pairs of points $(u,v)$ such that $u,v \in t\Newt$,
compute the integer point $pu-v$, and copy $\tau^{-1}$ of the term
$F_{pu-v} \bmod{p^{N}}$ from $F \bmod{p^{N}}$ into the correct
position in the matrix. Thus no arithmetic operations in the ring are
required here, except $W^{2}$ computations of $\tau^{-1} =
\tau^{a-1}$. These arithmetic operations can safely be ignored
since in Lemma \ref{Lem-ComputingF} we have already counted
$\Oh(|J|\tilde{W}^{2})$ ring operations.
Computation of the
appropriate indices $(u,v)$ does require one to find all lattice
points in certain polytopes and we return to that in Section
\ref{Sec-BitCompAux}.

Finally, for Step 3 we have the following estimate.

\begin{lemma}\label{Lem-ComputingApowers}
With the notation as in the statement of Theorem \ref{Thm-MDTF},
the product
\[ (\prod_{i = 0}^{\extdeg-1}\tau^{-i}(A_N))^{k} \]
can be computed given the matrix $A_N$ in
\[\Oh(W^{3}\log{(\extdeg k)})\]
operations in the ring $\Rp/(p^{N})$.
\end{lemma}

\begin{proof}
A fast square-and-multiply style algorithm may be used to compute
the power $\alpha^{\extdeg}$ in $\Oh(\log{\extdeg})$
``matrix ring operations''. Specifically,
working with the matrix representations one may
compute $\alpha^{r + s}$ from $\alpha^{r}$  and $\alpha^{s}$
using
\[ \prod_{i = 0}^{r + s - 1} \tau^{-i}(A_N) =
\prod_{i = 0}^{r - 1} \tau^{-i}(A_N) \left( \tau^{-r}(
\prod_{i = 0}^{s - 1} \tau^{-i}(A_N)) \right).\]
Now the case $r = s = 2^{c}$ for some $c$ gives us
the ``square'' step (computing $\alpha^{2^{c+1}}$ from $\alpha^{2^{c}}$) and
the case $r = 2^{c}$ and $s$ arbitrary the ``multiply'' step.
These two operations may be combined to give a
fast exponentiation method in a straightforward way.
The time required to compute a matrix for $\alpha^{\extdeg}$ from
one for $\alpha$ is thus $\Oh(\log{\extdeg})$ ``matrix ring
operations''. By matrix ring operations we mean
multiplication of matrices of size $W$ over
$\Rp/(p^{N})$, and also computing
$\tau^{-i}(B)$ for some $1 \leq i \leq \extdeg -1$
and matrix $B$ of this form.
The former requires $O(W^{3})$ operations in the ring
$\Rp/(p^{N})$ using standard algorithms.
Since $\tau^{-i} = \tau^{\extdeg -i}$ the latter may be computed
in $\Oh(W^{2})$ ring operations (that is, applications of a power of
$\tau$).
Thus the time for computing a matrix for
$\alpha^{\extdeg}$ from one for $\alpha$ is $\Oh(W^{3}\log{\extdeg})$.

Having obtained a matrix for
$\alpha^{\extdeg}$ one may then compute a matrix for $\alpha^{\extdeg
k}$ using
the standard square-and-multiply algorithm. This requires
$\Oh(\log{k})$ matrix multiplications, that is
$\Oh(W^{3}\log{k})$ operations in $\Rp/(p^{N})$.
Thus the total time required is as claimed.
\end{proof}

Note that the exponent $3$ for multiplication of matrices can be
improved to around $2.4$ using faster methods \cite[page 330]{GaGe}.

Gathering these results we find the following.

\begin{lemma}\label{Lem-RingOps}
The running time of Algorithm \ref{Alg-Torus} is
\[ \Oh(\tilde{W}^{2}|J| + W^{3}\log{(\extdeg k)})\]
ring operations. Here $|J|$ is the number of non-zero
terms in $f$, and $\tilde{W}$ and $W$ are defined as at
the start of Section \ref{Sec-RO}, with $N = \extdeg k(n+1)$.
\end{lemma}

\subsubsection{Bit complexity arising from ring operations}
\label{Sec-BitCompRing}

Using Lemma \ref{Lem-R0ops} one may now calculate the number of
bit operations in the algorithm which arise from operations
in the ring $\Rp/(p^{N})$. In this section we shall also
count the small contribution from computing the Teichm\"{u}ller
lifting of the coefficients of $f$, and also the
precomputation required for powers of $\tau$, which it was convenient
to ignore in Section \ref{Sec-RO},

\begin{definition}\label{Def-NormVol}
Let the polytope $\Newt$ from Definition \ref{Def-Deltas}
have dimension $\tilde{n} \leq n + 1$.
Let $V(\Newt)$ be the $\tilde{n}$-dimensional
volume of $\Newt$. Denote
by $\nvol = \tilde{n}!V(\Newt)$ the ``normalised'' volume of
$\Newt$.
\end{definition}

Since $f$ is non-zero we have $\tilde{n} \geq 1$ and certainly
$\nvol > 0$.

\begin{proposition}\label{Prop-BitComp}
The running time of Algorithm \ref{Alg-Torus} is
\[ \SoftOh(\extdeg^{3n+7}k^{3n+5}n^{3n+5}\nvol^{3}p^{2n+4})\]
bit operations, plus the contribution from operations
outside of the ring $\Rp/(p^{N})$.
The space complexity in bits is
\[ \SoftOh(\extdeg^{2n+4}k^{2n+3}n^{2n+3}\nvol^{2}p)\]
plus the contribution from operations outside of
$\Rp/(p^{N})$.
\end{proposition}

\begin{proof}
To compute the complexity first observe that $t = \lceil
p^{2}N/(p-1)^{2} \rceil - 1 \leq 4N$.
Thus $t,N = \Oh(\extdeg kn)$ in the algorithm. Also
$\tilde{t} = \Oh(pN) = \Oh(pt)$.
Now
\[
\begin{array}{rcl} W
& = & \#(t\Newt) \\
& \leq & \tilde{n}!V(t\Newt) + \tilde{n}\\
& = & \tilde{n}!t^{\tilde{n}}V(\Newt) + \tilde{n}\\
& = & \Oh(\nvol t^{n+1}).
\end{array}
\]
Here we use the Blichfeldt
bound $\#(P) \leq m!V(P) + m$
for any $m$-dimensional polytope $P$ from \cite[page 144]{GR},
the fact $V(t\Newt) = t^{\tilde{n}}V(\Newt)$
since $\Newt$ is $\tilde{n}$-dimensional, and also that
$\tilde{n} \leq n+1$.
Similarly $\tilde{W} = \#(\tilde{t}\Newt) =
\Oh(\nvol \tilde{t}^{n+1}) = \Oh(\nvol t^{n+1}p^{n+1})$.

Thus from Lemma \ref{Lem-RingOps} the number of ring operations is
\[ \Oh((\nvol t^{n+1}p^{n+1})^{2}(\nvol+n) +
(t^{n+1}\nvol)^{3}\log{(\extdeg k)})\]
since $|J| \leq \#(\Newt) - 1 \leq \nvol + n$.
Thus the bit complexity which arises
from ring operations is by Lemma \ref{Lem-R0ops}
\[ \Oh(((\nvol t^{n+1}p^{n+1})^{2}(\nvol + n) +
(t^{n+1}\nvol)^{3}\log{(\extdeg k)}
)(p\extdeg N\log{p})^{2}).\]
Tidying up and ignoring
logarithmic factors we get
\[ \SoftOh(\nvol^{3}t^{2n+2}p^{2n+4}\extdeg^{2}N^{2} +
\nvol^{3}t^{3n+3}p^{2}\extdeg^{2}
N^{2}).\]
The second term is dominant in all factors except $p$. For
simplicity we take the estimate of
\[ \SoftOh(\nvol^{3}t^{3n+3}p^{2n+4}\extdeg^{2}N^{2}).\]
Now putting $t,N = \Oh(\extdeg nk)$ we get
\[ \SoftOh(\nvol^{3}\extdeg^{3n+7}k^{3n+5}n^{3n+5}p^{2n+4}).\]
This is the total bit complexity which arises from ``ring
operations'', as defined at the start of Section \ref{Sec-RO}. There remains
the contribution from computing the Teichm\"{u}ller liftings of the
coefficients of $f$ in Step 1, and also the precomputation
required for the map $\tau$. By Lemma \ref{Lem-R0ops}, this is
easily seen to be absorbed in the above estimate.

With regard to the space complexity, this is dominated by the space
required to store the matrix $A_{N}$, which is $\Oh(W^{2})$ ring
elements. Putting $W = \Oh(\nvol(\extdeg nk)^{n+1})$ and using Lemma
\ref{Lem-R0ops} gives us the result.
\end{proof}

One may replace $n$ in the exponents in Proposition
\ref{Prop-BitComp} by $\tilde{n} - 1$;
however, this only gives an improvement when the Newton polytope
is not full-dimensional.

\subsubsection{Bit complexity arising from auxiliary operations}
\label{Sec-BitCompAux}

It remains to bound the complexity which arises from operations
outside of the ring $\Rp/(p^{N})$ in Steps 1 and 2.
In Step 1
manipulation of exponents of polynomials will add an extra term
$\Oh(\log(pNd))$ to the running time, which can safely be
ignored.

In Step 2 one is required to find all lattice points which
lie in $t\Newt$, for $t = \Oh(\extdeg nk)$. The complexity of
this step will depend upon the input polytope $\Newt$. For
a ``general'' $\Newt$ one requires methods from
computational convex geometry which are not in the spirit
of the present exposition. Thus our total bit complexity
estimate for Algorithm \ref{Alg-Torus} will just be as
in Proposition \ref{Prop-BitComp} ``plus the contribution
from finding all lattice points in $t\Newt$''.

At this stage for simplicity we shall restrict to
the choice of $\delta = \delta_{2}$ in Definition \ref{Def-Deltas}.
Thus we take $\Newt = \Newt_{2}$ as the convex hull in $\R^{n+1}$
of the origin and the $n+1$ points
\[ (1,0,\dots,0),(1,d,0,\dots,0),\dots,(1,0,\dots,0,d).\]
For this
case the required set of lattice points is
\[ \{ (r_{0},r_{1},\dots,r_{n})\,|\,r_{1} + \dots + r_{n} \leq
dr_{0} \leq dt\}\]
and so no computations are required here. Also, now
$\nvol$ equals $d^{n}$ and directly from Proposition \ref{Prop-BitComp}
we get the following result.

\begin{proposition}\label{Prop-BitCompTotal}
Let Algorithm \ref{Alg-Torus}$^{\prime}$ be exactly as Algorithm
\ref{Alg-Torus} only with input restricted to the choice of
polytope $\Newt = \Newt_{2}$ described in the preceding paragraph.
The total running time of Algorithm \ref{Alg-Torus}$^{\prime}$ is
\[ \SoftOh(\extdeg^{3n+7}k^{3n+5}n^{3n+5}d^{3n}p^{2n+4})\]
bit operations.
The total space complexity in bits is
\[ \SoftOh(\extdeg^{2n+4}k^{2n+3}n^{2n+3}d^{2n}p).\]
\end{proposition}

We shall use this restricted version of Algorithm \ref{Alg-Torus}
in the proofs of our main
results in the next section.

\subsection{Proofs of the results in the Introduction}\label{Sec-FinalProof}

To compute the number of points on the affine variety defined
by a polynomial $f$ one simply uses the torus decomposition
of $\fqk^{n}$. Specifically, for any subset $S \subseteq \{1,2,
\dots ,n\}$ let $G_{k}^{S}$ denote the set of points
\[ \{(x_{1},\dots,x_{n})\,|\,x_{i} \in \fqk, x_{i} = 0 \Leftrightarrow x \in
S\}.\]
Denote by $f^{S}$ the polynomial obtained from $f$ by setting
to zero all indeterminates $X_{i}$ which occur in $f$ for $i \in S$.
Denote by $N_{k}^{S}$ the number of solutions of $f^{S} = 0$
in the torus $G_{k}^{S}$ of dimension $n - |S|$. Then
$N_{k} = \sum_{S} N_{k}^{S}$ where the sum is over all subsets
of $\{1,2,\dots,n\}$. Each number $N_{k}^{S}$ can be computed
using Algorithm \ref{Alg-Torus}$^{\prime}$. (If some $f^{S}$ is identically
zero or has degree $0$ then $N_{k}^{S} = (q^{k} - 1)^{n - |S|}$ or
$0$, respectively, and Algorithm \ref{Alg-Torus}$^{\prime}$ is not required!)
Thus by $2^{n}$ applications of
this algorithm we obtain $N_{k}$ as desired.

Now to obtain the whole zeta function $Z(f/\fq)$
it suffices to count $N_{k}$ for all
$k = 1, \dots , \deg(r) + \deg(s)$, where
\[ Z(f/\fq)(T) = \frac{r(T)}{s(T)} \]
with $r$ and $s$ coprime polynomials in $1 + T\Z[T]$.
More precisely, it is enough to know upper bounds
$\deg(r) \leq D_1$ and $\deg(s) \leq D_2$, and compute
$N_{k}$ for $k = 1,\dots,D_1 + D_2$. Then
use the linear algebra method
described
prior to \cite[Corollary 2.8]{DW00}, which we now supplement
with further details.

Let $u(T) = 1 + \sum_{i = 1}^{D_{1}}u_{i}T^{i}$ and
$v(T) = 1 + \sum_{i = 1}^{D_{2}}v_{i}T^{i}$ have indeterminate
coefficients. Write $Z(f/\fq)(T) =
1 + z_{1}T + z_{2}T^{2} + \dots\,$. This power series
has non-negative integer coefficients and it can easily be computed
modulo $T^{D_{1} + D_{2} + 1}$ given $N_{k}$ for $k = 1,
\dots,D_{1} + D_{2}$.
The equation $v(T)Z(f/\fq) \equiv u(T) \bmod{T^{D_{1} + D_{2} + 1}}$ defines
a linear system $Ax = y$ where $A$ is a known
square $D_{1} + D_{2}$ integer matrix and $y$ a known
integer column vector. The entries in $A$ and $y$ are just coefficients
from the power series $Z(f/\fq) \bmod{T^{D_{1} + D_{2} + 1}}$. Let
$b$ be a bound on their bit length.
The unknown entries in $x$ are the coefficients of $u$ and $v$.
By \cite{DW00} the set of all solutions to this system consists of precisely
those vectors $x$ derived by specialising the coefficients
of $u(T)$ and $v(T)$ to equal those of
$d(T)r(T)$ and $d(T)s(T)$, respectively,
for some $d(T) \in 1 + T\Z[T]$ with degree at most $\min(D_{1}-\deg(r),
D_{2} - \deg(s))$. In particular, the system has a unique
solution (i.e. $\det(A) \ne 0$) if and only if either $\deg(r) = D_{1}$
or $\deg(s) = D_{2}$ (or both).
The determinant $\det(A)$ can be computed
using the small primes method in \cite[Algorithm 5.10]{GaGe} in a number
of bit operations bounded by $\SoftOh((D_{1} +
D_{2})^{4}b^{2})$ \cite[Theorem 5.12]{GaGe}.
Let us assume now that $\det(A) \ne 0$ and so the system
has a unique solution, namely the unknown
vector containing the integer coefficients of $r$ and $s$.
Let $B$ be a bound on the bit length of these coefficients.
Find the unique solution to the linear system modulo enough small
primes which do not divide $\det(A)$, and recover this integer solution
using the Chinese remainder theorem. Precisely, work modulo a collection of
such primes whose product has bit length greater than $B$.
This second step requires $\SoftOh((D_{1} + D_{2})^{4}B^{2})$ bit operations
using Gaussian elimination
(this can be improved with a Pad\'{e} Approximation algorithm
\cite[Section 5.9]{GaGe}).
Values for $b$ and $B$ may be deduced from the bound
$N_{k} \leq q^{nk}$. Specifically, one may show from this
that the absolute values of the
reciprocal zeros of $r$ and $s$ are all $\leq q^{n}$,
and so we can take $B = \Oh((D_{1} + D_{2})n\log{q})$. Also,
we can take $b = \Oh((D_{1} + D_{2})^{2}n\log{q})$.
If in the above we find $\det(A) = 0$ then we must have that
$\deg(r) < D_{1}$ and $\deg(s) < D_{2}$. In this case
one must first reduce $D_{2}$, say, and
compute determinants until the correct value $D_{2} = \deg(s)$ is
found (then $\det(A) \ne 0$ and the above method works).

By the refinement of Bombieri's
degree bound \cite{EB78}
from \cite[Equation (1.13)]{AS87}, the ``total degree''
$\deg(r) + \deg(s)$ is bounded by $2^{n+1}6^{n+1}(n+1)!V(\Newt_{1})$,
where $\Newt_{1}$ is the polytope in $\R^{n+1}$ derived from the
Newton polytope of $f$ (see the paragraph following
Definition \ref{Def-Deltas}).
Certainly $(n+1)!V(\Newt_{1}) \leq d^{n}$.
Hence we may take $D_{1},D_{2} = 2^{4n+4} d^{n}$ and
so $D_1 + D_2 = 2^{4n+5} d^{n}$.

\begin{theorem}\label{Thm-MainRefined}
Let $f$ be a polynomial in $n$ variables of total
degree $d > 0$ over $\fq$, where
$q = p^{\extdeg}$.
The full zeta function $Z(f/\fq)$ can be computed deterministically
in
\[ \SoftOh(2^{13n^{2}}\extdeg^{3n+7}d^{3n^{2}+9n}p^{2n+4}) \]
bit operations.
(Here we use soft-Oh notation which ignores
logarithmic factors, as defined at the start of Section
\ref{Sec-ComplexityAnalysis}.)
\end{theorem}

\begin{proof}
From Proposition \ref{Prop-BitCompTotal} and the torus decomposition
method, the bit complexity of computing $N_{k}$ for
$k = 1,\dots,2^{4n+5}d^{n}$ is
\[\SoftOh(\sum_{k = 1}^{2^{4n+5} d^{n}}
\{\extdeg^{3n+7}k^{3n+5}n^{3n+5}d^{3n}p^{2n+4}\}2^{n}).\]
(The contribution from recovering $Z(f/\fq)$ from the $N_{k}$ is
absorbed in this estimate.)
Tidying up the
factor in $n$ we get the
claimed result.
\end{proof}

Since we may assume that $d > 1$ we have $2^{n^{2}} = \Oh(d^{n^{2}})$,
and Theorem \ref{Thm-Main} now follows.

The proof of Corollary \ref{Cor-GeneralVar} was explained in the
introduction, and we finish with some comments on Corollary
\ref{Cor-Jacobian}.  By Weil's theorem, the zeta function of the
smooth projective curve $\tilde{V}$ from Corollary \ref{Cor-Jacobian} is
of the form
\[
Z(\tilde{V})(T)= \frac{P(T)}{(1-T)(1-qT)}
\]
for some polynomial $P(T)$ whose reciprocal roots have complex
absolute value $q^{1/2}$. Since the (possibly singular) affine
curve $V$ and the smooth projective curve $\tilde{V}$ differ in
only finitely many closed points, we deduce that the zeta function
of $V$ is of the form
\[
Z(V)(T)=\frac{P(T)Q(T)}{1-qT},
\]
where $Q(T)$ is a rational function whose zeros and poles are
roots of unity. This zeta function, and in particular the rational
function $P(T)Q(T)$, may be computed within the time bound in
Corollary \ref{Cor-Jacobian} by Corollary \ref{Cor-GeneralVar}. In
terms of the pure weight decomposition \cite{DW00}, the
polynomial $P(T)$ (respectively,
$Q(T)$) is exactly the pure weight $1$ (respectively, weight $0$) part of
the product $P(T)Q(T)$, and can be recovered quickly from
$P(T)Q(T)$ via the LLL polynomial factorization algorithm. In our
current special case, one can proceed directly without using the
LLL-factorization algorithm. By repeatedly removing the common
factor of the numerator of
$P(T)Q(T)$ with $T^{s} - 1$ for $\phi(s)$ (Euler totient
function) not greater than the total degree of $P(T)Q(T)$, the
desired polynomial $P(T)$ can be recovered. The order of the group
of rational points on $\tilde{V}$ is simply $P(1)$, see
\cite{DW00}. Thus for fixed dimension and finite field one can
compute the order of the group of rational points on the Jacobian
of a smooth projective curve in time polynomial in the degree $d$.

\end{document}